\input amstex
\documentstyle{amsppt}
\magnification=1200
\UseAMSsymbols
\font\tencursive=eurm10
\font\sevencursive=eurm7
\newfam\cursivefam
       \textfont\cursivefam=\tencursive
       \scriptfont\cursivefam=\sevencursive

\topmatter
\title
A uniform ergodic theorem for some N\"orlund means
\endtitle
\rightheadtext{A uniform ergodic theorem for some N\"orlund means}
\author Laura Burlando
\endauthor
\address
Dipartimento di Matematica dell'Universit\`a 
di Genova, Via Dodecaneso 35, 16146 Genova, ITALY
\endaddress
\email
burlando\@dima.unige.it
\endemail
\abstract
We obtain a uniform ergodic theorem for the sequence $\frac1{s(n)}
\sum_{k=0}^n(\varDelta s)(n-k)\,T^k$, where $\varDelta$ is the inverse 
of the endomorphism on the vector space of scalar sequences which maps 
each sequence into the sequence of its partial sums, $T$ is a bounded linear 
operator on a Banach space and $s$ is a divergent nondecreasing 
sequence of strictly positive real numbers, such that $\lim_{n\rightarrow+\infty}
s(n+1)/s(n)=1$ and $\varDelta^qs\in\ell_1$ for some positive integer 
$q$. Indeed, we prove that if $T^n/s(n)$ converges to zero in the uniform 
operator topology, then the sequence of averages above converges in 
the same topology if and only if 1 is either in the resolvent set of $T$, or a 
simple pole of the resolvent function of $T$.
\endabstract
\endtopmatter

\document
\heading{\bf1. Introduction}\endheading
Throughout this paper, \comment
when the scalar field is not specified, 
we assume it may be either $\Bbb C$ or $\Bbb R$ and denote it by 
$\Bbb K$. Also, 
\endcomment
we will write $\Bbb N$ and $\Bbb Z_+$ for the sets of 
nonnegative integers and of strictly positive integers, respectively. 
Also, for each $\nu\in\Bbb N$, we will wrire $\Bbb N_\nu$ for the set 
of all nonnegative integers $n$ satisfying $n\geq\nu$.
\par
$\Bbb K$ will stand for either $\Bbb R$ or $\Bbb C$, and we will denote 
by $\Bbb K^\Bbb N$ the vector space (over $\Bbb K$) of all sequences 
in $\Bbb K$. For each vector space $V$ over $\Bbb K$, let $0_V$ 
and $I_V$ denote respectively the zero element of $V$ and the identity 
operator on $V$. If $V$ and $W$ are vector spaces over $\Bbb K$ and 
$\varLambda:V\longrightarrow W$ is a linear map, let $\Cal N(\varLambda)$ 
and $\Cal R(\varLambda)$ stand respectively for the kernel and the 
range of $\varLambda$.
\par
For each normed space $X$, we will write ${\|\ \|}_X$ for the norm of 
$X$, and $L(X)$ for the normed algebra of all bounded linear 
operators on $X$. Henceforth, by {\sl convergence in} $L(X)$ of a 
sequence of bounded linear operators on $X$, we will mean 
convergence with respect to the topology induced by ${\|\ \|}_{L(X)}$, 
that is, the uniform operator topology.
\par
If $X$ is a complex nonzero Banach space, then $L(X)$ is a complex 
Banach algebra---with identity $I_X$. For each $T\in L(X)$, let 
$r(T)$ and $\sigma(T)$ stand respectively for the spectral radius and for 
the spectrum of $T$. Also, let $\rho(T)$ and $\goth R_T$ stand 
respectively for the resolvent set and for the resolvent function of $T$. 
Namely, $\rho(T)=\Bbb C\setminus\sigma(T)$ and
$\goth R_T:\rho(T)\ni\lambda\longmapsto(\lambda I_X-T)^{-1}\in
L(X)$. It is well known that $\goth R_T$ is analytic on the open 
set $\rho(T)$.
\par
In [D1], N. Dunford obtained several results about convergence 
of the sequence $f_n(T)$ in different topologies (where $T\in L(X)$ 
for a complex Banach space $X$, and, for each $n\in\Bbb N$, $f_n$ 
is a complex-valued function, holomorphic in some 
open neighborhood of $\sigma(T)$). The uniform ergodic 
theorem, establishing equivalence between convergence of the 
sequence $\frac1n\sum\limits_{k=0}^{n-1}T^k$ in $L(X)$ and $1$ 
being either in $\rho(T)$ or a simple pole of $\goth R_T$, under the 
hypothesis $\lim\limits_{n\rightarrow+\infty}\frac1n{\|T^n\|}_{L(X)}
=0$, is a special case of one of these results (see [D1], $3.16$; 
see also [D2], comments following Theorem $8$). Notice that if 
the sequence \,$\frac1n\sum\limits_{k=0}^{n-1}T^k$ converges in 
$L(X)$, then \,$\frac1n{\|T^n\|}_{L(X)}$ necessarily 
converges to zero, as \,$\frac1nT^n=\frac{n+1}n\Bigl(\frac1{n+1}\sum
\limits_{k=0}^nT^k\Bigr)-\frac1n\sum\limits_{k=0}^{n-1}T^k$ for 
each $n\in\Bbb Z_+$.
\par
More general means of the sequence of the 
iterates of the bounded linear operator $T$ than the arithmetical 
ones involved in the uniform ergodic theorem, that is, the $(C,\alpha)$ 
{\sl means} \,$\frac1{A_\alpha(n)}\sum\limits_{k=0}^nA_{\alpha-1}(n-k)\,T^k$, 
$n\in\Bbb N$ (where $\alpha\in(0,+\infty)$, and $A_\alpha$ and 
$A_{\alpha-1}$ denote respectively the sequences of Ces\`aro 
numbers---whose definition is recalled here in Section $2$---of order 
$\alpha$ and $\alpha-1$; notice that for $\alpha=1$ 
we have \,$\frac1{A_\alpha(n)}\sum\limits_{k=0}^nA_{\alpha-1}(n-k)\,
T^k=\frac1{n+1}\sum\limits_{k=0}^nT^k$ for each $n\in\Bbb N$), were 
considered by E. Hille in [Hi]. Indeed, in [Hi], Theorem $6$ he proved 
that if the sequence $\frac1{A_\alpha(n)}\sum\limits_
{k=0}^nA_{\alpha-1}(n-k)\,T^k$ converges to some $E\in L(X)$ in $L(X)$, then \,
$\frac{{\|T^n\|}_{L(X)}}{n^\alpha}\longrightarrow0$ as $n\rightarrow
+\infty$ and $\lim\limits_{\lambda\rightarrow1^+}{\|(\lambda-1)
\goth R_T(\lambda)-E\|}_{L(X)}=0$. Notice that the former of 
these two conditions yields $r(T)\leq1$, and then the latter is 
equivalent to $1$ being either in $\rho(T)$, or a simple pole of 
$\goth R_T$, and moreover $E$ being the residue of $\goth R_T$ 
at $1$ (see the result recorded here as Theorem $2.4$). 
Theorem $6$ of [Hi] also provides a partial converse of this, establishing that if 
$T$ is power-bounded and $\lim\limits_{\lambda\rightarrow1^+}
{\|(\lambda-1)\goth R_T(\lambda)-E\|}_{L(X)}=0$, then $\lim\limits_
{n\rightarrow+\infty}{\biggl\|\frac1{A_\alpha(n)}\sum\limits_
{k=0}^nA_{\alpha-1}(n-k)\,T^k-E\biggr\|}_{L(X)}=0$ for each $\alpha
\in(0,+\infty)$.
\par
More recently, an improvement of [Hi], Theorem $6$ was obtained by 
T. Yoshimoto, who in [Y], Theorem $1$ replaced power-boundedness 
of $T$ by $\lim\limits_{n\rightarrow+\infty}\frac{{\|T^n\|}_{L(X)}}
{n^\omega}=0$ (where $\omega=\min\{1,\alpha\}$). Finally, in [E], 
E. Ed-dari was able to complete the $(C,\alpha)$ uniform ergodic 
theorem, by proving that the sequence $\frac1{A_\alpha(n)}\sum\limits_
{k=0}^nA_{\alpha-1}(n-k)\,T^k$ converges to $E$ in $L(X)$ if 
and only if \,$\frac{{\|T^n\|}_{L(X)}}{n^\alpha}\longrightarrow0$ as 
$n\rightarrow+\infty$ and $\lim\limits_{\lambda\rightarrow1^+}
{\|(\lambda-1)\goth R_T(\lambda)-E\|}_{L(X)}=0$. E. Ed-dari's 
result is recorded here as Theorem $2.9$.
\par
We are interested here in obtaining a uniform ergodic theorem 
for the {\sl N\"orlund means} of the sequence $T^n$, that is, for 
the means \,$\frac1{s(n)}\sum\limits_{k=0}^n(\varDelta s)(n-k)\,T^k$, 
$n\in\Bbb N$, where $s$ is a divergent nondecreasing sequence of 
strictly positive real numbers (and $\varDelta:\Bbb K^\Bbb N
\rightarrow\Bbb K^\Bbb N$ is as in the abstract; see Definition 
$4.1$ here). Notice that for $s=A_\alpha$, $\alpha\in(0,+\infty)$, 
one obtains the $(C,\alpha)$ means.
\par
In Section $2$ we collect some preliminaries, in order to make 
this paper as self-contained as possible. In Sections $3$, $4$ and 
$5$ we derive some properties of real sequences, that we use in the 
final section dealing with bounded linear operators.
\par
\comment
Sections $3$ and $4$ contain some results---concerning real 
sequences---that we will use in the proof of our generalization of 
the uniform ergodic theorem. We suspect that some of these 
results---in particular, among the ones about the least concave 
majorant---may be known. However, we have not been able to find 
any bibliographical reference. Thus, also for the 
convenience of the reader, we include all the proofs here.
\endcomment
\newpage
In Section $3$ we are concerned with the least concave majorant 
of a real sequence. In particular, in Theorem $3.9$ we prove that 
if $b$ is a real sequence such that the sequence ${\bigl(\frac{b(n)}
n\bigr)}_{n\in\Bbb Z_+}$ is bounded from above and $b$ is not, 
then the least concave majorant of $b$, besides being strictly 
increasing and divergent, has a subsequence that is asymptotic to 
the corresponding subsequence of $b$.
\par
In Section $4$ we mainly deal with the real sequences $s$ for 
which $\varDelta^ps$ is concave for some $p\in\Bbb N$. The 
main result of this section is Theorem $4.7$, in which we derive 
several properties of a sequence $s$ of nonnegative real numbers 
such that $\varDelta^ps$ is concave and unbounded from above 
for some $p\in\Bbb N$. In particular, we prove that $s$ is strictly 
increasing and divergent, $\lim\limits_{n\rightarrow+\infty}\frac
{s(n+1)}{s(n)}=1$, and $\varDelta^{p+2}s\in\ell_1$. Also, in 
Example $4.9$ we show that if $\alpha\in(0,+\infty)$ the 
sequence $A_\alpha$ satisfies the hypotheses of Theorem $4.7$ 
(for $p=[\alpha]$ if $\alpha\notin\Bbb Z_+$; for $p=\alpha-1$ if 
$\alpha\in\Bbb Z_+$).
\par
In Section $5$ we introduce an index $\Cal H(b)$ ($\in\Bbb N
\cup\{+\infty\}$) for a real sequence $b$, such that $\Cal H(b)<+
\infty$ if and only if the sequence ${\bigl(\frac{b(n)}{n^m}\bigr)}_
{n\in\Bbb Z_+}$ is bounded from above for some $m\in\Bbb N$, 
in which case $\Cal H(b)$ is the minimum of such nonnegative 
integers $m$. In Theorem $5.3$ we use Theorem $3.9$ to 
prove that if $b$ is unbounded from above and such that 
$\Cal H(b)<+\infty$, then $b$ has a majorant $s$ which satisfies 
the hypotheses of Theorem $4.7$ for $p=\Cal H(b)-1$, and moreover 
is such that $\limsup\limits_{n\rightarrow+\infty}\frac{b(n)}{s(n)}\in
\bigl[\frac1{\Cal H(b)},1\bigr]$. We also prove (in Proposition $5.4$) 
that if $a$ is a real sequence such that $\varDelta^qa\in\ell_1$ 
for some $q\in\Bbb Z_+$, then $\Cal H(a)\leq q-1$.
\par
Section $6$ contains our main result, that is Theorem $6.7$: we 
prove that if $T$ is a bounded linear operator on a complex 
Banach space, and $b$ is a divergent sequence of strictly positive 
real numbers, such that $\Cal H(b)<+\infty$ and $\lim\limits_{n\rightarrow+\infty}\frac
{{\|T^n\|}_{L(X)}}{b(n)}=0$ (which gives $r(T)\leq1$), then, 
for each divergent nondecreasing 
sequence $s$ of strictly positive real numbers, such that $\lim\limits_
{n\rightarrow+\infty}\frac{s(n+1)}{s(n)}=1$, $\varDelta^qs\in\ell_1$ 
for some $q\in\Bbb N_2$, and the sequence ${\bigl(\frac{b(n)}
{s(n)}\bigr)}_{n\in\Bbb N}$ is bounded (which gives $\lim\limits_
{n\rightarrow+\infty}\frac{{\|T^n\|}_{L(X)}}{s(n)}=0$), the sequence 
$\frac1{s(n)}\sum\limits_{k=0}^n(\varDelta s)(n-k)\,T^k$ converges 
in $L(X)$ if and only if $1$ is either in $\rho(T)$, or a simple pole 
of $\goth R_T$. The sequence $s$ can be chosen so that it is 
not infinite of higher order than $b$, and $\varDelta^ps$ is 
concave and unbounded from above for some $p\in\Bbb N$. 
We conclude this section---and the paper---with an example (Example 
6.10), showing that, contrary to the case of the sequence $A_\alpha$ 
considered in Theorem $6$ of [Hi], convergence in $L(X)$ of 
the sequence $\frac1{s(n)}\sum\limits_{k=0}^n(\varDelta s)(n-k)\,T^k$ 
does not imply $\lim\limits_{n\rightarrow+\infty}\frac{{\|T^n\|}_{L(X)}}
{s(n)}=0$, even if $s$ satisfies the hypotheses of Theorem $4.7$.
\heading{\bf2. Preliminaries}\endheading
If $X$ is a Banach space, and $Y$, $Z$ are closed subspaces of $X$, 
satisfying $X=Y\oplus Z$, by the {\sl projection of $X$ onto $Y$ along} 
$Z$ we mean the bounded linear map $P:X\longrightarrow X$ such that 
$Px\in Y$ and $x-Px\in Z$ for every $x\in X$. Notice that $I_X-P$ is 
the projection of $X$ onto $Z$ along $Y$, and that $P^2=P$. On the 
other hand, if $E\in L(X)$ satisfies $E^2=E$, it is easily seen that 
$\Cal R(E)$ is closed in $X$, $X=\Cal R(E)\oplus\Cal N(E)$, and 
$E$ is the projection of $X$ onto $\Cal R(E)$ along $\Cal N(E)$.
\par
We begin by recalling a classical characterization of simple poles of 
$\goth R_T$, that will be useful to us in this paper.
\proclaim\nofrills{Theorem 2.1 \rm(see [TL], V, $10.1$, $10.2$, $6.2$, 
$6.3$ and $6.4$, and IV, $5.10$){\bf.}\usualspace}
Let $X$ be a complex nonzero Banach space, $T\in L(X)$ and 
$\lambda_0\in\Bbb C$. If $\lambda_0$ is a simple pole of $\goth
R_T$, then $\lambda_0$ is an eigenvalue of $T$, $\Cal
N\bigl({(\lambda_0I_X-T)}^n\bigr)=\Cal N(\lambda_0I_X-T)$ 
and $\Cal R\bigl({(\lambda_0I_X-T)}^n\bigr)=\Cal R(\lambda_0I_X-T)$ 
for every $n\in\Bbb Z_+$, $\Cal R(\lambda_0I_X-T)$ is closed in 
$X$, $X=\Cal N(\lambda_0I_X-T)\oplus\Cal R(\lambda_0I_X-T)$, 
and the projection of $X$ onto $\Cal N(\lambda_0I_X-T)$ along 
$\Cal R(\lambda_0I_X-T)$ coincides with the residue of $\goth R_T$ 
at $\lambda_0$. Conversely, if $X=\Cal N(\lambda_0I_X-T)\oplus\Cal
R(\lambda_0I_X-T)$, then $\lambda_0$ is either in $\rho(T)$, or 
else a simple pole of $\goth R_T$.
\endproclaim
If $X$ is a complex nonzero Banach space and $T\in L(X)$, following 
[TL] , Definition on page $310$, we denote by $\goth A(T)$ the set 
of all complex-valued holomorphic functions $f$ whose domain 
$\text{Dom}(f)$ is an open neighbourhood of $\sigma(T)$. For each 
$f\in\goth A(T)$, the operator $f(T)\in L(X)$ is defined as follows:
$$
f(T)=\frac1{2\pi i}\int_{+\partial D}f(\lambda)\,\goth R_T(\lambda)\,
d\lambda,
$$
where $+\partial D$ denotes the positively oriented boundary of $D$, 
and $D$ is any open bounded subset of $\Bbb C$, such that $D
\supseteq\sigma(T)$,  $\overline D\subseteq\text{Dom}(f)$, $D$ has 
a finite number of components, with pairwise disjoint closures, and 
$\partial D$ consists of a finite number of simple closed rectifiable curves, 
no two of which intersect; the integral above does not depend on the 
particular choice of $D$ (see [TL], comment $2$ on pages $310$--$311$; 
see also [D1], $2.2$, $2.3$ and $2.6$). We recall that for each polynomial 
$\goth p:\Bbb C\ni \lambda\longmapsto\sum\limits_{k=0}^na_k
\lambda^k\in\Bbb C$ (where $n\in\Bbb N$, and 
$a_0,\dots,a_n\in\Bbb C$), we have $\goth p(T)=\sum\limits_{k=0}^n
a_kT^k$ (see [TL], V, $8.1$).
\newline
We will use the following convergence result for the elements of 
$\goth A(T)$, due to N. Dunford, a special case of which is the classical 
uniform ergodic theorem.
\proclaim\nofrills{Theorem 2.2 \rm(see [D1], $3.16$){\bf.}\usualspace}
Let $X$ be a complex nonzero Banach space, $T\in L(X)$, and 
${(f_n)}_{n\in\Bbb N}$ be a sequence in $\goth A(T)$, satisfying 
$1\in\text{\rm Dom}(f_n)$ for each $n\in\Bbb N$, such that $\lim
\limits_{n\rightarrow+\infty}f_n(1)=1$ and $(I_X-T)f_n(T)\longrightarrow0_
{L(X)}$ in $L(X)$ as $n\rightarrow+\infty$. Then the following three conditions 
are equivalent:
\roster
\item"(2.2.1)"there exists $E\in L(X)$ such that $E^2=E$, $\Cal R(E)=\Cal
N(I_X-T)$, and $f_n(T)\longrightarrow E$ in $L(X)$ as $n\rightarrow
+\infty$;
\item"(2.2.2)"$1$ is either in $\rho(T)$, or a simple pole of $\goth R_T$;
\item"(2.2.3)"$\Cal R(I_X-T)$ is closed and $X=\Cal N(I_X-T)\oplus\Cal
R(I_X-T)$.
\endroster
\endproclaim
\remark{Remark 2.3}
We remark that, under the hypotheses of Theorem $2.2$, each of conditions 
$(2.2.1)$--$(2.2.3)$ is actually equivalent to each of the following two 
conditions (which at first glance might respectively appear to be weaker 
and stronger than them):
\roster
\item"(2.3.1)"{\it the sequence ${\bigl(f_n(T)\bigr)}_{n\in\Bbb N}$ converges 
in $L(X)$;}
\item"(2.3.2)" {\it$\Cal R(I_X-T)$ is closed, $X=\Cal N(I_X-T)\oplus\Cal
R(I_X-T)$, and the sequence ${\bigl(f_n(T)\bigr)}_{n\in\Bbb N}$ converges 
in $L(X)$ to the projection of $X$ onto $\Cal N(I_X-T)$ along $\Cal
R(I_X-T)$.}
\endroster
Equivalence between $(2.2.1)$ and $(2.3.1)$ is observed in [D2], 
comments following Theorem $8$. For the convenience of the reader, we 
give here a proof of equivalence of these five conditions. 
Indeed, it suffices to prove that $(2.3.1)$ implies $(2.3.2)$. Suppose that 
$(2.3.1)$ is satisfied, and let $E\in L(X)$ be such that 
$f_n(T)\longrightarrow E$ in $L(X)$ as $n\rightarrow+\infty$. We prove that 
then $E^2=E$ and $\Cal R(E)=\Cal N(I_X-T)$.
\newline
We begin by proving that for each $x\in\Cal N(I_X-T)$ we have $Ex=x$. 
This is clear if $\Cal N(I_X-T)=\{0_X\}$. If instead  $\Cal N(I_X-T)\neq\{0_X\}$, 
then $1\in\sigma(T)$, and $\goth R_T(\lambda)x=\frac1{\lambda-1}\,x$ 
for every $\lambda\in\rho(T)$. Hence (see [TL], V, $1.3$) $f_n(T)x=f_n(1)x$ 
for every $n\in\Bbb N$. Since $\lim\limits_{n\rightarrow+\infty}f_n(1)=1$, we 
conclude that $Ex=x$. This gives the desired result, which in turn yields \,
$\Cal N(I_X-T)\subseteq\Cal R(E)$. On the other hand, since $(I_X-T)E=
\lim\limits_{n\rightarrow+\infty}(I_X-T)f_n(T)=0_{L(X)}$, we have $\Cal R(E)
\subseteq\Cal N(I_X-T)$. Hence $\Cal R(E)=\Cal N(I_X-T)$, and $E^2=E$.
\newline
We have thus proved that the equivalent conditions $(2.2.1)$--$(2.2.3)$ 
are satisfied. Now we observe that, since $f_n(T)$ commutes with $I_X
-T$ for each $n\in\Bbb N$ by [TL], V, $8.1$, and consequently $E$ also 
does, we have $E(I_X-T)=0_{L(X)}$. Hence $\Cal R(I_X-T)\subseteq
\Cal N(E)$. Since $E^2=E$ and $\Cal R(E)=\Cal N(I_X-T)$ give 
$X=\Cal N(I_X-T)\oplus\Cal N(E)$, and condition $(2.2.3)$ in turn 
gives  $X=\Cal N(I_X-T)\oplus\Cal R(I_X-T)$, we conclude that $\Cal
N(E)=\Cal R(I_X-T)$. Then condition $(2.3.2)$ is satisfied.
\endremark
\vskip5pt
We also recall the following consequence of [D1], $3.16$.
\proclaim\nofrills{Theorem 2.4 \rm([E], $1.3$; [HP], $18.8.1$){\bf.}\usualspace}
Let $X$ be a complex nonzero Banach space and $T$, $E\in L(X)$. If 
there exists a sequence ${(\lambda_n)}_{n\in\Bbb N}$ in $\rho(T)$ 
such that $\lim\limits_{n\rightarrow+\infty}\lambda_n=1$ and \,
$(\lambda_n-1)\goth R_T(\lambda_n)\longrightarrow E$ in $L(X)$ 
as $n\rightarrow+\infty$, then $1$ is either in $\rho(T)$, or a simple 
pole of $\goth R_T$. Furthermore, $\Cal R(I_X-T)$ is closed in $X$, 
$X=\Cal N(I_X-T)\oplus\Cal  R(I_X-T)$ and $E$ is the projection 
of $X$ onto $\Cal N(I_X-T)$ along $\Cal R(I_X-T)$.
\endproclaim
For each $\alpha\in\Bbb R$, let $A_\alpha:\Bbb N\rightarrow\Bbb R$ 
denote the sequence of the {\sl Ces\`aro numbers of order} $\alpha$. 
That is,
$$
A_\alpha(n)=\binom{n+\alpha}n=\cases1&\text{if }n=0\\\frac{\prod\limits_
{j=1}^n(\alpha+j)}{n!}&\text{if }n\in\Bbb Z_+.\endcases
$$
Hence $A_\alpha(n)>0$ for each $n\in\Bbb N$ if $\alpha>-1$. Notice 
also that $A_0(n)=1$ for all $n\in\Bbb N$. We recall that
\roster
\item"(2.5)"$\sum\limits_{k=0}^nA_\alpha(k)=A_{\alpha+1}(n)$ \,{\it for 
each $n\in\Bbb N$ and each} $\alpha\in\Bbb R$ 
\endroster
and
\roster
\item"(2.6)"$\lim\limits_{n\rightarrow+\infty}\frac{A_\alpha(n)}{n^\alpha}
=\frac1{\Gamma(\alpha+1)}$ \ {\it for each} $\alpha\in\Bbb R\setminus
\{-k:k\in\Bbb Z_+\}$,
\endroster
where $\Gamma$ denotes Euler's gamma function (see for instance 
[Z], III, (1-11) and (1-15)).
\newline
The following well known identity---which we will need in the sequel---can 
be obtained from $(2.5)$ as a straightforward consequence, or else is 
not difficult to check directly, by induction on $n$.
\roster
\item"(2.7)"$\sum\limits_{k=j}^n\binom kj=\binom{n+1}{j+1}$ {\it \,for every 
$j\in\Bbb N$ and every} $n\in\Bbb N_j$.
\endroster
\remark{Remark 2.8}
Let $X$ be a complex nonzero Banach space, and let $T\in L(X)$. We 
recall that if the sequence ${\Bigl(\frac{{\|T^n\|}_{L(X)}}{n^\alpha}\Bigr)}_
{\!n\in\Bbb Z_+}$ is bounded for some $\alpha\in(0,+\infty)$, then $r(T)\leq1$. 
Indeed, if $M\in(0,+\infty)$ is such that \,$\frac{{\|T^n\|}_{L(X)}}{n^\alpha}\leq
M$ for each $n\in\Bbb Z_+$, then
$$
r(T)=\lim_{n\rightarrow+\infty}{{\|T^n\|}_{L(X)}^{\frac1n}}=\lim_{n\rightarrow
+\infty}{\Biggl(\frac{{\|T^n\|}_{L(X)}}{n^\alpha}\Biggr)}^{\!\!\frac1n}\leq
\lim_{n\rightarrow+\infty}M^{\frac1n}=1.
$$
\endremark
\vskip5pt
Finally, by also taking Theorem $2.4$ into account, the improvement of 
E. Hille's $(C,\alpha)$ ergodic theorem obtained by 
E. Ed-dari can be formulated as follows.
\proclaim\nofrills{Theorem 2.9 \rm(see [E], Theorem $1$){\bf.}\usualspace}
Let $X$ be a complex nonzero Banach space, $T\in L(X)$, and 
$\alpha\in(0,+\infty)$. Then, given any $E\in L(X)$, we have
$$
\lim_{n\rightarrow+\infty}{\left\|\,\frac{\sum\limits_{k=0}^nA_{\alpha-1}(n-k)\,T^k}
{A_\alpha(n)}-E\,\right\|}_{L(X)}=0
$$
if and only if
$$
\lim_{n\rightarrow+\infty}\frac{{\|T^n\|}_{L(X)}}{n^\alpha}=0\qquad\text{and}
\qquad\lim_{\lambda\rightarrow1^+}{\|(\lambda-1)\goth R_T(\lambda)-E\|}_
{L(X)}=0.\text{\footnotemark}
$$
\footnotetext{We point out that, by virtue of Remark $2.8$, $\lim\limits_
{n\rightarrow+\infty}\frac{{\|T^n\|}_{L(X)}}{n^\alpha}=0$ \ gives \ $r(T)\leq1$. 
Then $\rho(T)$ contains all real numbers $\lambda$ satisfying $\lambda>1$, 
which allows the limit $\lim\limits_{\lambda\rightarrow1^+}{\|(\lambda-1)\goth
R_T(\lambda)-E\|}_{L(X)}$ to be considered.}
Hence the following two conditions are equivalent:
\roster
\item"(2.9.1)"the sequence \ ${\left(\frac{\sum\limits_{k=0}^nA_{\alpha-1}
(n-k)\,T^k}{A_\alpha(n)}\right)}_{\!\!n\in\Bbb N}$ \ converges in $L(X)$;
\item"(2.9.2)"$\lim\limits_{n\rightarrow+\infty}\frac{{\|T^n\|}_{L(X)}}{n^\alpha}
=0$ and $1$ is either in $\rho(T)$, or a simple pole of $\goth R_T$.
\endroster
\endproclaim
\heading{\bf3. The least concave majorant of a real sequence}
\endheading
We begin with some results concerning the least concave majorant of a 
real sequence.
\par
We recall that a real sequence $a:\Bbb N\rightarrow\Bbb R$ is called 
{\sl concave} ({\sl convex}) if the real sequence ${\bigl(a(n+1)-a(n)\bigr)}_
{n\in\Bbb N}$ is nonincreasing (nondecreasing). Notice that $a$ is 
concave (convex) if and only if $a(n+1)\geq\frac{a(n)+a(n+2)}2$ 
($a(n+1)\leq\frac{a(n)+a(n+2)}2$\,) for every $n\in\Bbb N$.
\definition{Definition 3.1}
For each real sequence $a:\Bbb N\rightarrow\Bbb R$, let  
$\phi_a:[0,+\infty)\rightarrow\Bbb R$ be the function defined by
$$
\phi_a(x)=a(n)+(x-n)\bigl(a(n+1)-a(n)\bigr)\qquad\text{for every }x\in
[n,n+1]\text{ and every }n\in\Bbb N.
$$
\enddefinition
Notice that $\phi_a(x)=a(n)(n+1-x)+a(n+1)(x-n)$ for every $x\in[n,n+1]$ 
and every $n\in\Bbb N$. Hence $\phi_a(n)=a(n)$ for every $n\in\Bbb N$.
\proclaim{Proposition 3.2}
Let $a:\Bbb N\rightarrow\Bbb R$ be a real sequence. Then $a$ is 
concave if and only if the function $\phi_a$ is concave.
\endproclaim
\demo{Proof}
It is easily seen that $a$ is concave if $\phi_a$ is. Conversely, 
suppose $a$ to be concave. Notice that $\phi_a$ is continuous. Also, 
the right derivative $(\phi_a)_+'$ of $\phi_a$ exists at every point of 
$[0,+\infty)$, and $(\phi_a)_+'(x)=a(n+1)-a(n)$ for every $x\in[n,n+1)$ and 
every $n\in\Bbb N$. Since $a$ is concave, it follows that $(\phi_a)_+'$ 
is nonincreasing, and consequently (see [R], 5, Proposition 18) $\phi_a$ 
is concave.
\qed
\enddemo
We recall that a {\sl majorant} of a real sequence $b:\Bbb N\rightarrow\Bbb
R$ is a real sequence $c:\Bbb N\rightarrow\Bbb R$ satisfying $c(n)\geq
b(n)$ for every $n\in\Bbb N$.
\newline
The following result is probably known. Indeed, for instance, the authors 
of [AR] seem to be aware of it when (in the proof of Proposition $2.1$) 
they derive that the sequence ${(\rho_n)}_{n\in\Bbb N}$ has a least 
concave majorant from being $\lim\limits_{n\rightarrow+\infty}\frac{\rho_n}
n=0$. Anyway, we give a (short) proof here, for the convenience of the reader.
\proclaim{Proposition 3.3}
A real sequence $b:\Bbb N\rightarrow\Bbb R$ has a concave majorant 
if and only if the sequence ${\bigl(\frac{b(n)}n\bigr)}_{n\in\Bbb Z_+}$ 
is bounded from above.
\endproclaim
\demo{Proof}
By virtue of Proposition $3.2$, it is easily seen that $b$ has a concave 
majorant if and only if there exists a concave function $f:[0,+\infty)
\rightarrow\Bbb R$ such that $f(x)\geq\phi_b(x)$ for every $x\in
[0,+\infty)$. The latter condition is satisfied if and only if there 
exist $\alpha$, $\beta\in\Bbb R$ such that $\phi_b(x)\leq\alpha+
\beta x$ for every $x\in[0,+\infty)$ (see [G], Theorem $1.2$) or, equivalently, 
$b(n)\leq\alpha+\beta n$ for every $n\in\Bbb N$. Now it is straightforward 
to observe that such $\alpha$ and $\beta$ exist if and only if the 
sequence ${\bigl(\frac{b(n)}n\bigr)}_{n\in\Bbb Z_+}$ 
is bounded from above.
\comment
Suppose the sequence ${\bigl(\frac{b(n)}n\bigr)}_{n\in\Bbb Z_+}$ 
to be bounded from above and let $M\in\Bbb R$ be such that \ 
$\frac{b(n)-b(0)}n\leq M$ for all $n\in\Bbb Z_+$. Then $b(n)\leq b(0)+Mn$ 
for every $n\in\Bbb N$. If $a:\Bbb N\rightarrow\Bbb R$ is defined by
$$
a(n)=b(0)+Mn\qquad\text{for each }n\in\Bbb N,
$$
it follows that $a$ is a majorant of $b$ and, since
$$
a(n+1)-a(n)=M(n+1)-Mn=M\qquad\text{for every }n\in\Bbb N,
$$
$a$ is concave.
\newline
Conversely, suppose that $b$ has a concave majorant. Then there 
exists a concave sequence $a:\Bbb N\rightarrow\Bbb R$ such that 
$a(n)\geq b(n)$ for every $n\in\Bbb N$. Consequently, for each 
$n\in\Bbb Z_+$, we have
$$
\frac{b(n)}n\leq\frac{a(n)}n=\frac{a(0)+\sum\limits_{k=0}^{n-1}\bigl(a(k+1)-
a(k)\bigr)}n.\tag{3.3.1}
$$
Since the sequence ${\bigl(a(n+1)-a(n)\bigr)}_{n\in\Bbb N}$ is 
nonincreasing, it follows that there exists $\ell\in[-\infty,+\infty)$ such 
that $\lim\limits_{n\rightarrow+\infty}\bigl(a(n+1)-a(n)\bigr)=\ell$. 
Since \ $\frac{a(0)}n\longrightarrow0$ as $n\rightarrow+\infty$, from
the classical Ces\`aro means theorem and from $(3.3.1)$ we conclude 
that
$$
\ell=\lim_{n\rightarrow+\infty}\frac{\sum\limits_{k=0}^{n-1}
\bigl(a(k+1)-a(k)\bigr)}n=\lim_{n\rightarrow+\infty}\frac{a(0)+\sum
\limits_{k=0}^{n-1}\bigl(a(k+1)-a(k)\bigr)}n=\lim_{n\rightarrow+\infty}
\frac{a(n)}n.
$$
Since $\ell\in[-\infty,+\infty)$, it follows that the sequence ${\bigl(\frac
{a(n)}n\bigr)}_{n\in\Bbb Z_+}$ is bounded from above. By virtue of 
$(3.3.1)$, we conclude that also the sequence ${\bigl(\frac
{b(n)}n\bigr)}_{n\in\Bbb Z_+}$ is bounded from above.
\endcomment
We have thus obtained the desired result.
\qed
\enddemo
\remark{Remark 3.4}
If $a:\Bbb N\rightarrow\Bbb R$ is a concave sequence, then the 
sequence ${\bigl(\frac{a(n)}n\bigr)}_{n\in\Bbb Z_+}$ is bounded from 
above.
\endremark
\remark{Remark 3.5}
If a real sequence $b:\Bbb N\rightarrow\Bbb R$ has a concave 
majorant, then $b$ has a {\sl least} concave majorant $c$. Furthermore, 
we have
$$
c(n)=\inf\{a(n):a\in\Bbb R^{\Bbb N},\ a\text{ concave majorant of }b\}
\qquad\text{for every }n\in\Bbb N.
$$
Indeed, once one observes that the real sequence $c$ defined 
as above is concave and is a majorant of $b$, from the definition 
of $c$ it follows that each concave majorant of $b$ is also a 
majorant of $c$, that is, $c$ is the least concave 
majorant of $b$.
\endremark
\proclaim{Theorem 3.6}
Let $b:\Bbb N\rightarrow\Bbb R$ be a real sequence such that the 
sequence ${\bigl(\frac{b(n)}n\bigr)}_{n\in\Bbb Z_+}$ is bounded 
from above, and let $c:\Bbb N\rightarrow\Bbb R$ be the least concave 
majorant of $b$. Then $c$ satisfies the following properties.
\roster
\item"(3.6.1)"$c(0)=b(0)$.
\item"(3.6.2)"$c(n+1)=c(n)+\sup\left\{\frac{b(k)-c(n)}{k-n}:k\in\Bbb N_
{n+1}\right\}$ for every\footnote{Notice that, since the 
sequence ${\bigl(\frac{b(n)}n\bigr)}_{n\in\Bbb Z_+}$ is bounded from above, 
this supremum is finite for every $n\in\Bbb N$.} $n\in\Bbb N$.
\item"(3.6.3)"For each $n\in\Bbb N$ and each $k\in\Bbb N_{n+2}$, we have 
$$
\frac{b(k)-c(n+1)}{k-n-1}\leq\frac{b(k)-c(n)}{k-n}.
$$
If in addition
$$
\frac{b(k)-c(n)}{k-n}=\max\left\{\frac{b(j)-c(n)}{j-n}:j\in\Bbb N_{n+1}\right\},
$$
then
$$
\align
\frac{b(k)-c(h)}{k-h}&=\max\left\{\frac{b(j)-c(h)}{j-h}:j\in\Bbb N_{h+1}\right\}
\\&=\frac{b(k)-c(n)}{k-n}\qquad\text{for all \,}h=n,\dots,k-1,
\endalign
$$
and consequently \,$c(h+1)-c(h)=c(n+1)-c(n)$ for all \,$h=n,\dots,k-1$.
\endroster
\endproclaim
\demo{Proof}
Let $a:\Bbb N\rightarrow\Bbb R$ be defined by
$$
\multline
a(0)=b(0),\\
a(n+1)=a(n)+\sup\left\{\frac{b(k)-a(n)}{k-n}:k\in\Bbb N_{n+1}\right\}
\qquad\text{for every }n\in\Bbb N.
\endmultline
$$
We prove that $a=c$. First of all, we prove that $a(n)\geq b(n)$ for every 
$n\in\Bbb N$.
\newline
We proceed by induction. The desired result clearly holds for 
$n=0$. Besides, if for some $n\in\Bbb N$ we have $a(n)\geq b(n)$, then 
$$
a(n+1)=a(n)+\sup\left\{\frac{b(k)-a(n)}{k-n}:k\in\Bbb N_{n+1}\right\}
\geq a(n)+b(n+1)-a(n)=b(n+1).
$$
We have thus proved that $a$ is a majorant of $b$. Now we prove that 
$a$ is concave.
\newline
For each $n\in\Bbb N$, we have
$$
a(n+1)-a(n)=\sup\left\{\frac{b(k)-a(n)}{k-n}:k\in\Bbb N_{n+1}\right\}
\tag{3.6.4}
$$
and
$$
a(n+2)-a(n+1)=\sup\left\{\frac{b(k)-a(n+1)}{k-n-1}:k\in\Bbb N_{n+2}
\right\}.\tag{3.6.5}
$$
Now, for each $k\in\Bbb N_{n+ 2}$, $(3.6.4)$ yields
$$
\multline
\hskip10pt\frac{b(k)-a(n+1)}{k-n-1}=\frac{b(k)-a(n)}{k-n-1}-\frac{a(n+1)-a(n)}{k-n-1}\\
\hskip54pt=\left(\frac{k-n}{k-n-1}\right)\left(\frac{b(k)-a(n)}{k-n}\right)-\left(\frac
1{k-n-1}\right)\sup\left\{\frac{b(j)-a(n)}{j-n}:j\in\Bbb N_{n+1}\right\}
\\\leq\left(\frac{k-n}{k-n-1}\right)\left(\frac{b(k)-a(n)}{k-n}\right)-
\frac1{k-n-1}\left(\frac{b(k)-a(n)}{k-n}\right)\\=\left(\frac{b(k)-a(n)}{k-n}\right)
\left(\frac{k-n}{k-n-1}-\frac1{k-n-1}\right)=\frac{b(k)-a(n)}{k-n}.
\endmultline
\tag3.6.6
$$
Now from $(3.6.6)$, together with $(3.6.4)$ and $(3.6.5)$, we derive that 
$a$ is concave. Hence $a$ is a concave majorant of $b$. In order to 
conclude that $a=c$, it suffices to prove that $c(n)\geq a(n)$ for every 
$n\in\Bbb N$. We proceed by induction.
\comment
Let $x:\Bbb N\rightarrow\Bbb R$ be a concave majorant of $b$. Proceeding 
by induction on $n$, we prove that $x(n)\geq c(n)$ for every $n\in\Bbb N$.
\endcomment
\newline
Since $c(0)\geq b(0)=a(0)$, the desired inequality holds for $n=0$. Now let 
$n\in\Bbb N$ be such that $c(n)\geq a(n)$. Since $c$ is concave and 
is a majorant of $b$, from Proposition $3.2$ and from the three chord lemma 
we conclude that for each $k\in\Bbb N_{n+1}$ we have
\comment
$$
b(k)-x(n)\leq x(k)-x(n)=\sum_{j=n}^{k-1}\bigl(x(j+1)-x(j)\bigr)\leq(k-n)
\bigl(x(n+1)-x(n)\bigr)
$$
\endcomment
$$
c(n+1)-c(n)\geq\frac{c(k)-c(n)}{k-n}\geq\frac{b(k)-c(n)}{k-n}
$$
and consequently, since $c(n)\geq a(n)$,
$$
\multline
c(n+1)\geq c(n)+\frac{b(k)-c(n)}{k-n}=c(n)+\frac{b(k)-a(n)}{k-n}+\frac
{a(n)-c(n)}{k-n}\\=\frac{(k-n-1)}{k-n}\,c(n)+\frac1{k-n}\,a(n)+\frac{b(k)-a(n)}
{k-n}\\\geq\left(\frac{k-n-1}{k-n}+\frac1{k-n}\right)a(n)+\frac{b(k)-a(n)}
{k-n}=a(n)+\frac{b(k)-a(n)}{k-n}.
\endmultline
$$
Then
$$
c(n+1)\geq a(n)+\sup\left\{\frac{b(k)-a(n)}{k-n}:k\in\Bbb N_{n+1}\right\}
=a(n+1).
$$
We have thus proved that $c(n)\geq a(n)$ for every $n\in\Bbb N$. Then 
$a=c$, from which we obtain $(3.6.1)$ and $(3.6.2)$. Also, the 
inequality in $(3.6.3)$ now follows from $(3.6.6)$.
\newline
In order to complete the proof, it remains to prove that if $n\in\Bbb N$ 
and $k\in\Bbb N_{n+2}$ satisfy $\frac{b(k)-c(n)}{k-n}=
\max\left\{\frac{b(j)-c(n)}{j-n}:j\in\Bbb N_{n+1}\right\}$, then
$\frac{b(k)-c(h)}{k-h}=\max\left\{\frac{b(j)-c(h)}{j-h}:j\in\Bbb N_{h+1}\right\}
=\frac{b(k)-c(n)}{k-n}$ \ for all $h=n,\dots,k-1$.
\newline
Let $n\in\Bbb N$, $k\in\Bbb N_{n+2}$ be as above. As a straightforward 
consequence of $(3.6.6)$, we obtain
$$
\frac{b(k)-c(n+1)}{k-n-1}=\frac{b(k)-c(n)}{k-n}
$$
and consequently
$$
c(n+2)-c(n+1)=\sup\left\{\frac{b(j)-c(n+1)}{j-n-1}:j\in\Bbb N_{n+2}\right\}
\geq\frac{b(k)-c(n)}{k-n}=c(n+1)-c(n).
$$
Since $c$ is concave, the opposite inequality also holds. Hence
$$
\frac{b(k)-c(n+1)}{k-n-1}=\max\left\{\frac{b(j)-c(n+1)}{j-n-1}:j\in\Bbb
N_{n+2}\right\}=\frac{b(k)-c(n)}{k-n}.
$$
If $k=n+2$, the proof is complete. Otherwise, we finish the proof 
by applying the same argument again \,$k-n-2$ times, with \,
$\frac{b(k)-c(n)}{k-n}$ \,replaced by \,$\frac{b(k)-c(h)}{k-h}$, 
$h=n+1,\dots,k-2$.
\newline
\qed
\enddemo
\proclaim{Lemma 3.7}
Let $a:\Bbb N\rightarrow\Bbb R$ be a real sequence. If there exists 
$\nu\in\Bbb N$ such that $a(\nu)\geq\limsup\limits_{n\rightarrow
+\infty}a(n)$, then the set $\{a(n):n\in\Bbb N\}$ has a maximum.
\endproclaim
\demo{Proof}
We set $\ell=\limsup\limits_{n\rightarrow+\infty}a(n)$ and observe that 
$\ell\in[-\infty,+\infty)$. Since assuming $a(n)\leq\ell$ for every $n\in\Bbb
N$ yields $\ell\in\Bbb R$ and $a(\nu)=\ell\geq a(n)$ for every $n\in\Bbb
N$---which means that $a(\nu)$ is the maximum of $\{a(n):n\in\Bbb
N\}$, we may assume that $a(\nu)>\ell$. Then there exists $n_0
\in\Bbb N$ such that $a(n)<a(\nu)$ for every $n\in\Bbb N_{n_0}$, 
from which we conclude that $\nu<n_0$. Now let $n_1\in\{0,\dots,
n_0-1\}$ be such that $a(n_1)\geq a(k)$ for all $k=0,\dots,n_0-1$. 
It suffices to remark that for each $n\in\Bbb N_{n_0}$ we have 
$a(n)<a(\nu)\leq a(n_1)$. Hence $a(n_1)\geq a(n)$ for every 
$n\in\Bbb N$.
\qed
\enddemo
\proclaim{Theorem 3.8}
Let $b:\Bbb N\rightarrow\Bbb R$ be a real sequence such that the 
sequence ${\bigl(\frac{b(n)}n\bigr)}_{n\in\Bbb Z_+}$ is bounded from above, 
and let $c:\Bbb N\rightarrow\Bbb R$ be the least concave majorant of $b$.
\roster
\item"(3.8.1)"If we set $\ell=\limsup\limits_{n\rightarrow+\infty}\,\frac{b(n)}n$, 
we have $\ell\in[-\infty,+\infty)$, \ $c(n+1)-c(n)\geq\ell$ \ for every $n\in
\Bbb N$ and\ \,$\limsup\limits_{k\rightarrow+\infty}\left(\frac{b(k)-c(n)}{k-n}\right)
=\ell$ for every $n\in\Bbb N$.
\item"(3.8.2)"If we set
$$
\goth N=\{0\}\,\cup\,\biggr\{n\in\Bbb Z_+:c(n)-c(n-1)=
\max\left\{\tsize\frac{b(k)-c(n-1)}{k-n+1}:k\in\Bbb N_n\right\}\biggr\},
$$
it follows that \ $n\in\goth N$ $\Longrightarrow$ $\{0,\dots,n\}\subseteq\goth N$.
\item"(3.8.3)"If we set $N=\sup(\goth N)$ and ${(\nu_k)}_{k\in\Bbb N}$ is the 
nondecreasing sequence of nonnegative integers defined by \ $\nu_0=0$,
\item""$\nu_{k+1}=\cases\min\left\{n\in\Bbb N_{\nu_k+1}:
\frac{b(n)-c(\nu_k)}{n-\nu_k}=c(\nu_k+1)-c(\nu_k)\right\}&\text
{if\footnotemark\ }\nu_k<N\\\nu_k&\text{if }\nu_k\geq N\endcases$
\footnotetext{Notice that if $\nu_k<N$, from $(3.8.2)$ it follows that 
$\nu_k+1\in\goth N$ and consequently $c(\nu_k+1)-c(\nu_k)=
\frac{b(n)-c(\nu_k)}{n-\nu_k}$ for some $n\in\Bbb N_{\nu_k+1}$.}
\newline
for every $k\in\Bbb N$, it follows that $\nu_k\in\goth N$ for every $k\in\Bbb N$ 
and $\{\nu_k:k\in\Bbb N\}=\mathbreak\{n\in\Bbb N:c(n)=b(n)\}$.
\item"(3.8.4)"If $\goth N$ is finite (that is, $N\in\Bbb N$, $N=\max(\goth
N)$), then the sequence ${(\nu_k)}_{k\in\Bbb N}$ is eventually 
constant (and consequently $\nu_k\geq N$ for sufficiently large $n$). 
Furthermore, if we set $k_0=min\{k\in\Bbb N:\nu_k\geq N\}$ we have: 
$\nu_k=N$ for every $k\in\Bbb N_{k_0}$, $\nu_k<\nu_{k+1}$ for each 
$k\in\Bbb N$ satisfying $k<k_0$, $\ell\in\Bbb R$, and, for each $n\in\Bbb N$,
\item""$c(n)=\cases b(\nu_k)+\left(\frac{n-\nu_k}{\nu_{k+1}-\nu_k}\right)
\bigl(b(\nu_{k+1})-b(\nu_k)\bigr)&\text{if }\nu_k\leq n\leq\nu_{k+1}\text{ for 
some }k\in\Bbb N\\&\hskip71.59pt\text{satisfying }k<k_0\\b(N)+\ell(n-N)&\text{if }n
\geq N
\endcases$
\item""$=\cases\left(\frac{\nu_{k+1}-n}{\nu_{k+1}-\nu_k}\right)
b(\nu_k)+\left(\frac{n-\nu_k}{\nu_{k+1}-\nu_k}\right)b(\nu_{k+1})
&\text{if }\nu_k\leq n\leq\nu_{k+1}\text{ for some }k\in\Bbb N\\&\hskip88.23pt
\text{satisfying }k<k_0\\b(N)+\ell(n-N)&\text{if }n\geq N.
\endcases$
\item""Finally, $c(n)>b(n)$ for every $n\in\Bbb N_{N+1}$.
\item"(3.8.5)"If $\goth N$ is infinite (that is, $\goth N=\Bbb N$), then the 
sequence ${(\nu_k)}_{k\in\Bbb N}$ is strictly increasing. Furthermore, 
for each $k\in\Bbb N$ we have
$$
\align
c(n)&=b(\nu_k)+\left(\frac{n-\nu_k}{\nu_{k+1}-\nu_k}\right)
\bigl(b(\nu_{k+1})-b(\nu_k)\bigr)\\&=\left(\frac{\nu_{k+1}-n}{\nu_{k+1}-\nu_k}\right)
b(\nu_k)+\left(\frac{n-\nu_k}{\nu_{k+1}-\nu_k}\right)b(\nu_{k+1})
\endalign
$$
for every $n\in\Bbb N$ satisfying $\nu_k\leq n\leq\nu_{k+1}$.
\endroster
\endproclaim
\demo{Proof}
We begin by proving $(3.8.1)$. As a straightforward consequence of 
${\bigl(\frac{b(n)}n\bigr)}_{n\in\Bbb Z_+}$ being bounded from above, we 
have $\ell\in[-\infty,+\infty)$.
\newline
Now let ${(k_j)}_{j\in\Bbb N}$ be a strictly 
increasing sequence of strictly positive integers such that 
$\lim\limits_{j\rightarrow+\infty}\frac{b(k_j)}{k_j}=\ell$. Then for each $n\in\Bbb
N$ we have
$$
\lim_{j\rightarrow+\infty}\left(\frac{b(k_j)-c(n)}{k_j-n}\right)=\lim_{j\rightarrow+\infty}
\left(\frac{k_j}{k_j-n}\right)\left(\frac{b(k_j)}{k_j}-\frac{c(n)}{k_j}\right)=\ell.
\tag3.8.6
$$
Hence, by virtue of $(3.6.2)$,
$$
c(n+1)-c(n)=\sup\left\{\frac{b(k)-c(n)}{k-n}:k\in\Bbb N_{n+1}\right\}\geq\ell.
$$
Now if we set $s=\limsup\limits_{k\rightarrow+\infty}\left(\frac{b(k)-c(n)}
{k-n}\right)$, from $(3.8.6)$ we also derive that $s\geq\ell$. On the other 
hand, if ${(m_j)}_{j\in\Bbb N}$ is a strictly increasing sequence of 
nonnegative integers such that $\lim\limits_{j\rightarrow+\infty}
\left(\frac{b(m_j)-c(n)}{m_j-n}\right)=s$, we obtain
$$
\frac{b(m_j)}{mj}=\left(\frac{b(m_j)-c(n)}{m_j-n}\right)\left(\frac{m_j-n}
{m_j}\right)+\frac{c(n)}{m_j}@>>j\rightarrow+\infty>s,
$$
which gives $s\leq\ell$. Hence $s=\ell$.
\newline
$(3.8.1)$ is thus proved. Now we prove $(3.8.2)$.
\newline
Fix $n\in\goth N$ and let $k\in\{0,\dots,n\}$. We prove that $k\in\goth N$. This 
is clearly true if $k=0$. If $k\in\Bbb Z_+$, then $n\in\Bbb Z_+$ and \ $c(n)-c(n-1)=
\frac{b(p_n)-c(n-1)}{p_n-n+1}$ \ for some $p_n\in\Bbb N_n$.
It is not restrictive to assume $k\leq n-1$ (which gives \ $n-k\in\Bbb Z_+$, 
$n\in\Bbb N_2$, $k-1\leq n-2$). Since for each $j\in\{k-1,\dots,n-2\}$ 
we have $j+2\leq p_n$, and consequently $\frac{b(p_n)-c(j+1)}
{p_n-j-1}\leq\frac{b(p_n)-c(j)}{p_n-j}$ by $(3.6.3)$, by taking $(3.8.1)$ 
into account we obtain
$$
\frac{b(p_n)-c(k-1)}{p_n-k+1}\geq\frac{b(p_n)-c(n-1)}{p_n-n+1}=
c(n)-c(n-1)\geq\ell=\limsup\limits_{m\rightarrow+\infty}\left(\frac{b(m)-c(k-1)}
{m-k+1}\right).
$$
Now from Lemma $3.7$ we conclude that the set $\left\{\frac{b(m)-c(k-1)}
{m-k+1}:m\in\Bbb N_k\right\}$ has a maximum. This, together with $(3.6.2)$, 
yields $k\in\goth N$.
\newline
We prove $(3.8.3)$.
\newline
We begin by proving that for each $k\in\Bbb N$ we have $\nu_k\in\goth
N$ and $\bigl\{n\in\{0,\dots\nu_k\}:c(n)=b(n)\bigr\}\mathbreak=\{\nu_j:j=0,\dots,k\}$. 
We proceed by induction. We set
$$
\Cal S=\Bigl\{k\in\Bbb N:\nu_k\in\goth N\text{ and }
\bigl\{n\in\{0,\dots,\nu_k\}:c(n)=b(n)\bigr\}=\{\nu_j:j=0,\dots,k\}\Bigr\}.
$$
Since $\nu_0=0$, by the definition of $\goth N$ and by $(3.6.1)$ 
we have $0\in\Cal S$. Now suppose $k\in\Cal S$. Then $\nu_k\in\goth N$ 
and $\bigl\{n\in\{0,\dots\nu_k\}:c(n)=b(n)\bigr\}=\{\nu_j:j=0,\dots,k\}$. If 
$\nu_k\geq N$, we have $\nu_{k+1}=\nu_k\in\goth N$ and $\bigl\{n\in
\{0,\dots,\nu_{k+1}\}:c(n)=b(n)\bigr\}=\bigl\{n\in\{0,\dots\nu_k\}:c(n)=b(n)\bigr\}
=\{\nu_j:j=0,\dots,k\}=\{\nu_j:j=0,\dots,k+1\}$, which gives $k+1\in\Cal S$. 
Thus, let us assume $\nu_k<N$. Then $\nu_{k+1}\geq\nu_k+1$ and
$$
\nu_{k+1}=\min\left\{n\in\Bbb N_{\nu_k+1}:\frac{b(n)-c(\nu_k)}{n-\nu_k}=
c(\nu_k+1)-c(\nu_k)\right\}.\tag3.8.7
$$
From $(3.6.3)$ we derive that for each $j\in\Bbb N$ satisfying $\nu_k\leq
j\leq\nu_{k+1}-1$ we have
$$
\frac{b(\nu_{k+1})-c(j)}{\nu_{k+1}-j}=\max\left\{\frac{b(m)-c(j)}{m-j}:m\in
\Bbb N_{j+1}\right\}=\frac{b(\nu_{k+1})-c(\nu_k)}{\nu_{k+1}
-\nu_k}.\tag3.8.8
$$
By letting $j=\nu_{k+1}-1$, from $(3.8.8)$---together with $(3.6.2)$---we 
conclude that $\nu_{k+1}\in\goth N$ and besides
$$
b(\nu_{k+1})-c(\nu_{k+1}-1)=\max\left\{\frac{b(m)-c(\nu_{k+1}-1)}
{m-\nu_{k+1}+1}:m\in\Bbb N_{\nu_{k+1}}\right\}=
c(\nu_{k+1})-c(\nu_{k+1}-1),
$$
which gives $c(\nu_{k+1})=b(\nu_{k+1})$. Finally, for each $n\in\Bbb
N_{\nu_k+1}$ satisfying $n<\nu_{k+1}$, $(3.6.2)$, $(3.6.3)$, $(3.8.7)$ 
and $(3.8.8)$ give
$$
\multline
c(n)-c(n-1)=\max\left\{\frac{b(m)-c(n-1)}{m-n+1}:m\in\Bbb N_n\right\}\\=
\frac{b(\nu_{k+1})-c(\nu_k)}{\nu_{k+1}-\nu_k}>\frac{b(n)-c(\nu_k)}{n-\nu_k}
\geq b(n)-c(n-1),
\endmultline
$$
the latter inequality being trivially an equality if $n=\nu_k+1$, and being 
a consequence of $(3.6.3)$ if $n\geq\nu_k+2$ (as \,$\frac{b(n)-c(j)}
{n-j}\geq\frac{b(n)-c(j+1)}{n-j-1}$ for all $j=\nu_k,\dots,n-2$). Consequently, 
$c(n)>b(n)$. Hence
$$
\multline
\bigl\{n\in\{0,\dots,\nu_{k+1}\}:c(n)=b(n)\bigr\}\\=\bigl\{n\in\{0,\dots\nu_k\}:
c(n)=b(n)\bigr\}\cup\bigl\{n\in\{\nu_k+1,\dots,\nu_{k+1}\}:c(n)=b(n)\bigr\}
\\=\{\nu_j:j=0,\dots,k\}\cup\{\nu_{k+1}\}=\{\nu_j:j=0,\dots,k+1\},
\endmultline
$$
which gives $k+1\in\Cal S$.
\newline
We have thus proved that $\nu_k\in\goth N$ for every $k\in\Bbb N$. Also,
$$
\bigl\{n\in\{0,\dots,\nu_k\}: c(n)=b(n)\bigr\}=\{\nu_j: j=0,\dots,k\}\qquad
\text{for every }k\in\Bbb N.\tag3.8.9
$$
Now we prove that $\{n\in\Bbb N: c(n)=b(n)\}=\{\nu_k:k\in\Bbb N\}$.
\newline
If $N=+\infty$, then $\nu_{k+1}\geq\nu_k+1$ for every $k\in\Bbb N$, 
which gives $\lim\limits_{k\rightarrow+\infty}\nu_k=+\infty$, and 
consequently $\bigcup\limits_{k\in\Bbb N}\{0,\dots,\nu_k\}=\Bbb N$. 
Hence
$$
\multline
\{n\in\Bbb N:c(n)=b(n)\}=\bigcup_{k\in\Bbb N}\bigl\{n\in\{0,\dots,\nu_k\}:
c(n)=b(n)\bigr\}\\=\bigcup_{k\in\Bbb N}\{\nu_j:j=0,\dots,k\}=\{\nu_k:k\in
\Bbb N\},
\endmultline
$$
which is the desired result.
\newline
If $N<+\infty$, then there exists $\overline k\in\Bbb N$ such that 
$\nu_{\overline k}\geq N$: otherwise, if $\nu_k<N$ for all $k\in\Bbb
N$, the sequence ${(\nu_k)}_{k\in\Bbb N}$ would be strictly 
increasing and consequently we would have $+\infty=\lim\limits_
{k\rightarrow+\infty}\nu_k\leq N$, a contradiction. Hence $\nu_k
=\nu_{\overline k}$ for every $k\in\Bbb N_{\overline k}$. 
Furthermore, for each $n\in\Bbb N_{\nu_{\overline k}+1}$, we 
have $n>N$ and consequently $n\notin\goth N$. Then $n\in
\Bbb Z_+$ and, by virtue of $(3.6.2)$, $c(n)-c(n-1)>b(n)-c(n-1)$, 
which gives $c(n)>b(n)$. From this, together with $(3.8.9)$, 
we obtain
$$
\multline
\{n\in\Bbb N:c(n)=b(n)\}=\bigl\{n\in\{0,\dots,\nu_{\overline k}\}:c(n)
=b(n)\bigr\}\\=\{\nu_j:j=0,\dots,\overline k\,\}=\{\nu_k:k\in\Bbb N\}.
\endmultline
$$
We have thus finished the proof of $(3.8.3)$.
\newline
We prove $(3.8.4)$. Suppose $\goth N$ to be finite. Then $N\in
\Bbb N$ and $N=\max(\goth N)$. Also, we have already 
observed---in the proof of $(3.8.3)$---that the sequence ${(\nu_k)}
_{k\in\Bbb N}$ is eventually constant and not less than $N$. 
We set $k_0=\min\{k\in\Bbb N:\nu_k\geq N\}$. Then $\nu_{k_0}
\geq N$. On the other hand, since $\nu_{k_0}\in\goth N$ by 
$(3.8.3)$, we have $\nu_{k_0}\leq N$. Then $\nu_{k_0}=N$, 
and consequently $\nu_k=N$ for each $k\in\Bbb N_{k_0}$. 
Furthermore, for each $k\in\Bbb N$ satisfying $k<k_0$ we 
have $\nu_k<N$, and consequently $\nu_k<\nu_{k+1}$. From 
$(3.8.8)$ and $(3.6.2)$ we conclude that for each $n\in\Bbb N$ 
satisfying $\nu_k+1\leq n\leq\nu_{k+1}$ we have
$$
c(j)-c(j-1)=\frac{b(\nu_{k+1})-c(\nu_k)}{\nu_{k+1}-\nu_k}\qquad
\text{for all }j=\nu_k+1,\dots,n
$$
and consequently
$$
c(n)-c(\nu_k)=\sum_{j=\nu_k+1}^n\bigl(c(j)-c(j-1)\bigr)=(n-\nu_k)
\left(\frac{b(\nu_{k+1})-c(\nu_k)}{\nu_{k+1}-\nu_k}\right).
$$
Hence
$$
c(n)-c(\nu_k)=(n-\nu_k)\left(\frac{b(\nu_{k+1})-c(\nu_k)}
{\nu_{k+1}-\nu_k}\right)\qquad\text{for all }n=\nu_k,\dots,\nu_
{k+1}.\tag3.8.10
$$
Since $c(\nu_k)=b(\nu_k)$ by $(3.8.3)$, from $(3.8.10)$ 
we derive that for each $n\in\{\nu_k,\dots,\nu_{k+1}\}$ we 
have
$$
c(n)=c(\nu_k)+(n-\nu_k)\left(\frac{b(\nu_{k+1})-c(\nu_k)}
{\nu_{k+1}-\nu_k}\right)=b(\nu_k)+\left(\frac{n-\nu_k}{\nu_{k+1}-\nu_k}\right)
\bigl(b(\nu_{k+1})-b(\nu_k)\bigr).
$$
Now we prove that $\ell\in\Bbb R$ and $c(n)=b(N)+\ell(n-N)$ 
for every $n\in\Bbb N_N$.
\newline
For each $n\in\Bbb N_N$, we have $n+1>N$ and consequently 
$n+1\notin\goth N$. From $(3.6.2)$ we conclude that the set 
$\left\{\frac{b(k)-c(n)}{k-n}:
k\in\Bbb N_{n+1}\right\}$ has no maximum. From Lemma $3.7$ 
and from $(3.8.1)$ we derive that \,$\frac{b(k)-c(n)}{k-n}<\ell$ 
for every $k\in\Bbb N_{n+1}$, and consequently $\ell\in\Bbb R$. 
Besides, from $(3.6.2)$ we obtain
$$
c(n+1)-c(n)=\sup\left\{\frac{b(k)-c(n)}{k-n}:k\in\Bbb N_{n+1}\right\}
\leq\ell,
$$
which,  together with $(3.8.1)$, gives $c(n+1)-c(n)=\ell$. 
Notice also that by virtue of $(3.8.3)$, $\nu_{k_0}=N$ yields 
$c(N)=c(\nu_{k_0})=b(\nu_{k_0})=b(N)$. Now, proceeding by 
induction, we conclude that $c(n)=b(N)+\ell(n-N)$ for every 
$n\in\Bbb N_N$. Finally, from $(3.8.3)$ we derive that $c(n)
>b(n)$ for every $n\in\Bbb N_{N+1}$ and the proof of $(3.8.4)$ 
is complete.
\newline
We prove $(3.8.5)$. If we assume $\goth N$ to be infinite 
(or equivalently, by virtue of $(3.8.2)$, $\goth N=\Bbb N$), 
then $N=+\infty$ and consequently the sequence ${(\nu_k)}
_{k\in\Bbb N}$ is strictly increasing. The remaining assertion 
can be derived from $(3.8.8)$, $(3.6.2)$ and $(3.8.3)$, 
proceeding as in the proof of $(3.8.4)$. The proof is now finished.
\qed
\enddemo
\proclaim{Theorem 3.9}Let $b:\Bbb N\rightarrow\Bbb R$ be 
a real sequence such that the sequence ${\bigl(\frac{b(n)}n\bigr)}_
{n\in\Bbb Z_+}$ is bounded from above and $b$ is not, and let  
$c$ be the least concave majorant of $b$. Then $c$ is strictly 
increasing, $\lim\limits_{n\rightarrow+\infty}c(n)=+\infty$ and 
$\limsup\limits_{n\rightarrow+\infty}\frac{b(n)}{c(n)}=1$.
\endproclaim
\demo{Proof}
We set $\ell=\limsup\limits_{n\rightarrow+\infty}\frac{b(n)}n$. 
We observe that $\limsup\limits_{n\rightarrow+\infty}b(n)=+\infty$, 
and consequently $\ell\in[0,+\infty)$. From $(3.8.1)$ it follows 
that $c$ is nondecreasing, and consequently there exists 
$\lim\limits_{n\rightarrow+\infty}c(n)$. Since $c(n)\geq b(n)$ 
for every $n\in\Bbb N$, we conclude that
$$
\lim_{n\rightarrow+\infty}c(n)\geq\limsup_{n\rightarrow+\infty}b(n)
=+\infty.
$$
Hence $c(n)\longrightarrow+\infty$ as $n\rightarrow+\infty$. 
Now we prove that $c$ is strictly increasing. If $c$ were not 
strictly increasing, then---being $c$ nondecreasing---there 
would be $n_0\in\Bbb N$ such that $c(n_0+1)-c(n_0)=0$. 
Since $c$ is concave as well as nondecreasing, we would 
conclude that $c$ is eventually constant, in contradiction 
with $\lim\limits_{n\rightarrow+\infty}c(n)=+\infty$.
\newline
Finally, we prove that $\limsup\limits_{n\rightarrow+\infty}
\frac{b(n)}{c(n)}=1$. By virtue of Theorem $3.8$, one of the 
following two conditions is satisfied:
\roster
\item"(3.9.1)"there exists a strictly increasing sequence 
${(\nu_k)}_{k\in\Bbb N}$ of nonnegative integers such that 
$c(\nu_k)=b(\nu_k)$ for every $k\in\Bbb N$;
\item"(3.9.2)"there exists $N\in\Bbb N$ such that $c(n)>b(n)$ 
for every $n\in\Bbb N_{N+1}$ and $c(n)=b(N)+\ell(n-N)$ for
every $n\in\Bbb N_N$.
\endroster
If $(3.9.1)$ holds, it suffices to observe that $\limsup\limits_
{n\rightarrow+\infty}\frac{b(n)}{c(n)}\geq\lim\limits_
{k\rightarrow+\infty}\frac{b(\nu_k)}{c(\nu_k)}=1$. The opposite 
inequality follows from $c$ being a majorant of $b$.
\newline
If $(3.9.2)$ holds, then $\lim\limits_{n\rightarrow+\infty}
\bigl(b(N)+\ell(n-N)\bigr)=\lim\limits_{n\rightarrow+\infty}
c(n)=+\infty$ gives $\ell\in(0,+\infty)$. Hence
$$
\limsup_{n\rightarrow+\infty}\left(\frac{b(n)}{c(n)}\right)=
\limsup_{n\rightarrow+\infty}\left(\frac{b(n)}n\right)\cdot
\frac1{\ell+\left(\frac{b(N)-\ell N}n\right)}=1.
$$
The desired result is thus proved.
\qed
\enddemo
The following is a consequence of Remark $3.4$ and Theorem $3.9$. 
Alternatively, it can be derived from Proposition $3.2$ and the properties 
of concave functions.
\proclaim{Corollary 3.10}
If a concave sequence $a:\Bbb N\rightarrow\Bbb R$ is not bounded 
from above, then $a$ is strictly increasing and $\lim\limits_{n\rightarrow
+\infty}a(n)=+\infty$.
\endproclaim
\heading{\bf4. Real sequences with concave $p^{\text{th}}$-difference}
\endheading
\definition{Definition 4.1}
Let $\varSigma,\varDelta:\Bbb K^\Bbb N\rightarrow\Bbb K^\Bbb N$ 
be the linear operators defined by
$$
(\varSigma a)(n)=\sum_{k=0}^na(k)\qquad\text{and}\qquad(\varDelta 
a)(n)=\cases a(0)&\text{if }n=0\\a(n)-a(n-1)&\text{if }n\in\Bbb Z_+
\endcases
$$
for every $n\in\Bbb N$ and every $a\in\Bbb K^\Bbb N$.
\enddefinition
Notice that both linear operators $\varSigma$ and $\varDelta$ are bijective. 
Besides, $\varDelta=\varSigma^{-1}$ (or, equiva\-lently, $\varSigma=
\varDelta^{-1}$). We also remark that $\varDelta(\ell_1)\subseteq\ell_1$. 
Finally, we observe that the operator $\varSigma$ preserves inequalities: 
indeed, if $a$, $b\in\Bbb R^\Bbb N$ satisfy $a(n)\leq b(n)$ for each 
$n\in\Bbb N$, then $(\varSigma a)(n)\leq(\varSigma b)(n)$ for each $n
\in\Bbb N$.
\par
The following is a consequence of Proposition $3.2$ and of the the three 
chord lemma.
\proclaim{Lemma 4.2}
Let $a:\Bbb N\rightarrow\Bbb R$ be a concave sequence. Then
$$
n\bigl(a(k)-a(0)\bigr)\geq k\bigl(a(n)-a(0)\bigr)\qquad\text{for every }
n\in\Bbb N\text{ and every }k=0,\dots,n.
$$
\endproclaim
\proclaim{Theorem 4.3}
Let $a:\Bbb N\rightarrow\Bbb R$ be a nondecreasing concave sequence. 
Then for each $p\in\Bbb N$ we have
$$
\frac1{p+1}\binom{n+p}pa(n)+\frac p{p+1}\binom{n+p}pa(0)\leq
(\varSigma^p a)(n)\leq\binom{n+p}pa(n)\qquad\text{for every }n\in\Bbb N.
$$
\endproclaim
\demo{Proof}
We begin by proving that $(\varSigma^pa)(n)\leq\binom{n+p}pa(n)$ for 
all $n$, $p\in\Bbb N$. We proceed by induction on $p$.
\newline
For $p=0$, the desired inequality trivially holds for every $n\in\Bbb N$. Now 
let $p\in\Bbb N$ be such that $(\varSigma^pa)(n)\leq\binom{n+p}pa(n)$ for 
every $n\in\Bbb N$. Then for each $n\in\Bbb N$, since $a$ is nondecreasing 
we have
$$
(\varSigma^{p+1}a)(n)=\sum_{k=0}^n(\varSigma^pa)(k)\leq
\sum_{k=0}^n{\tsize\binom{k+p}p}a(k)\leq a(n)\sum_{k=0}^n\tsize{\binom{k+p}p}
=\binom{n+p+1}{p+1}a(n)
$$
by $(2.7)$.
\newline
We have thus proved the desired inequality. Now, proceeding again by 
induction on $p$, we prove that $(\varSigma^pa)(n)\geq\frac1{p+1}
\binom{n+p}pa(n)+\frac p{p+1}\binom{n+p}pa(0)$ for all $n$, $p\in\Bbb
N$.
\newline
For $p=0$, the desired inequality is trivially satisfied for every $n\in\Bbb N$. 
Now let $p\in\Bbb N$ be such that $(\varSigma^pa)(n)\geq\frac1{p+1}
\binom{n+p}pa(n)+\frac p{p+1}\binom{n+p}pa(0)$ for every $n\in\Bbb N$. 
We prove that $(\varSigma^{p+1}a)(n)\geq\frac1{p+2}\binom{n+p+1}{p+1}
a(n)+\frac{p+1}{p+2}\binom{n+p+1}{p+1}a(0)$ for every $n\in\Bbb N$.
\newline
For $n=0$, since $(\varSigma^{p+1}a)(0)=a(0)$ the desired inequality 
is straightforward. Now fix $n\in\Bbb Z_+$. Then from Lemma $4.2$ 
and $(2.7)$ we obtain
$$
\multline
(\varSigma^{p+1}a)(n)=\sum_{k=0}^n(\varSigma^pa)(k)\geq\frac1{p+1}
\sum_{k=0}^n\binom{k+p}pa(k)+\frac{p}{p+1}\,a(0)\sum_{k=0}^n
\binom{k+p}p\\
\allowdisplaybreak\geq\frac1{n(p+1)}\sum_{k=0}^n\binom{k+p}p
\bigl(ka(n)+(n-k)a(0)\bigr)+\frac p{p+1}\,a(0)\binom{n+p+1}{p+1}\\
\allowdisplaybreak
=\frac{\bigl(a(n)-a(0)\bigr)}n\sum_{k=1}^n\frac k{p+1}\binom{k+p}p
+\frac{a(0)}{p+1}\sum_{k=0}^n\binom{k+p}p+\frac p{p+1}\,a(0)
\binom{n+p+1}{p+1}\\
\allowdisplaybreak=\frac{\bigl(a(n)-a(0)\bigr)}n\sum_{k=1}^n
\binom{k+p}{p+1}+a(0)\binom{n+p+1}{p+1}\\
\allowdisplaybreak=\frac{\bigl(a(n)-a(0)\bigr)}n
\sum_{k=0}^{n-1}\binom{k+p+1}{p+1}+a(0)\binom{n+p+1}{p+1}\\
\allowdisplaybreak
=\frac{\bigl(a(n)-a(0)\bigr)}n\binom{n+p+1}{p+2}+a(0)\binom{n+p+1}
{p+1}\\\allowdisplaybreak=
\frac{\bigl(a(n)-a(0)\bigr)}{p+2}\cdot\frac{(n+p+1)!}{n\,(p+1)!\,(n-1)!}
+a(0)\binom{n+p+1}{p+1}\\\allowdisplaybreak
=\frac{\bigl(a(n)-a(0)\bigr)}{p+2}\binom{n+p+1}
{p+1}+a(0)\binom{n+p+1}{p+1}\\\allowdisplaybreak
=\frac1{p+2}\binom{n+p+1}{p+1}a(n)+
\frac{p+1}{p+2}\binom{n+p+1}{p+1}a(0),
\endmultline
$$
which is the desired result. The proof is now complete.
\qed
\enddemo
\proclaim{Lemma 4.4}
Let $b:\Bbb N\rightarrow\Bbb R$ be a nondecreasing sequence, 
satisfying $b(0)\geq0$ and $b(1)>0$. Then for each $k\in\Bbb Z_+$ 
the sequence $\varSigma^kb$ is convex and strictly increasing.
\endproclaim
\demo{Proof}
It suffices to prove that $\varSigma b$ is convex and strictly increasing. 
Indeed, once this is proved, the desired result 
follows by induction on $k$ (as $(\varSigma^kb)(0)=b(0)\geq0$ for 
every $k\in\Bbb Z_+$, and then $\varSigma^kb$ strictly increasing gives 
$(\varSigma^kb)(1)>0$).
\newline
Sice $b$ is nondecreasing, $b(1)>0$ yields $b(n)>0$ for every $n\in
\Bbb Z_+$, and consequently $(\varSigma b)(n+1)=(\varSigma b)(n)+
b(n+1)>(\varSigma b)(n)$ for every $n\in\Bbb N$. Hence $\varSigma
b$ is strictly increasing. Furthermore, being the sequence 
${\bigl((\varSigma b)(n+1)-(\varSigma b)(n)\bigr)}_{n\in\Bbb N}=
{\bigl(b(n+1)\bigr)}_{n\in\Bbb N}$ nondecreasing, $\varSigma b$ is 
convex. We have thus obtained the desired result.
\qed
\enddemo
\proclaim{Lemma 4.5}
Let $c:\Bbb N\rightarrow\Bbb R$ be a concave nondecreasing sequence. 
Then $\varDelta c$ is convergent, and $\varDelta^2c\in\ell_1$.
\endproclaim
\demo{Proof}
Since $c$ is concave and nondecreasing, it follows that the sequence 
${\bigl((\varDelta c)(n+1)\bigr)}_{n\in\Bbb N}$ is nonincreasing and 
$(\varDelta c)(n)\geq0$ for each $n\in\Bbb Z_+$. Then $\lim\limits_
{n\rightarrow+\infty}(\varDelta c)(n)=\lambda$ for some $\lambda\in
[0,+\infty)$ (and so $\varDelta c$ converges). Besides, $(\varDelta^2c)(n)
\leq0$ for each $n\in\Bbb N_2$. Since
$$
\sum_{k=0}^n(\varDelta^2c)(n)=(\varSigma\varDelta^2c)(n)=
(\varDelta c)(n)@>>n\rightarrow+\infty>\lambda
$$
(that is, the series $\sum\limits_{n=0}^{+\infty}(\varDelta^2c)(n)$ converges), 
being $\varDelta^2c$ eventually nonpositive 
it follows that the series $\sum\limits_{n=0}^{+\infty}|(\varDelta^2c)(n)|$ 
also converges. Hence $\varDelta^2c\in\ell_1$.
\qed
\enddemo
\remark{Remark 4.6}
Let $s\in\Bbb K^\Bbb N$ be an eventually nonzero sequence. If we fix 
$\nu\in\Bbb Z_+$ such that $s(n)\neq0$ for all $n\in\Bbb N_\nu$, then 
for each $n\in\Bbb N_\nu$ we have \,$\frac{(\varDelta s)(n)}{s(n)}=
1-\frac{s(n-1)}{s(n)}$. Hence
$$
\lim_{n\rightarrow+\infty}\frac{s(n+1)}{s(n)}=1\qquad\Longleftrightarrow
\qquad\lim_{n\rightarrow+\infty}\frac{(\varDelta s)(n)}{s(n)}=0.
$$
\endremark
\proclaim{Theorem 4.7}
Let $s:\Bbb N\rightarrow\Bbb R$ be a real sequence satisfying $s(0)
\geq0$, and let $p\in\Bbb N$ be such that the sequence $\varDelta^ps$ 
is concave and is not bounded from above. Then:
\roster
\item"(4.7.1)"$s(n)>0$ for every $n\in\Bbb Z_+$;
\item"(4.7.2)"$s$ is strictly increasing;
\item"(4.7.3)"$\lim\limits_{n\rightarrow+\infty}\frac{s(n)}{n^p}=+\infty$ 
(and consequently $\lim\limits_{n\rightarrow+\infty}s(n)=+\infty$);
\item"(4.7.4)"the sequence ${\bigl(\frac{s(n)}{n^{p+1}}\bigr)}_{n\in\Bbb
Z_+}$ is bounded;
\item"(4.7.5)"$\lim\limits_{n\rightarrow+\infty}\frac{s(n+1)}{s(n)}=1$;
\item"(4.7.6)"$\varDelta^{p+2}s\in\ell_1$.
\endroster
If in addition $p\in\Bbb Z_+$, then the sequence $s$ is also convex.
\endproclaim
\demo{Proof}
From Corollary $3.10$ it follows that  $\varDelta^ps$ is strictly 
increasing and $\lim\limits_{n\rightarrow+\infty}(\varDelta^ps)(n)
=+\infty$. Also, $(\varDelta^ps)(0)=s(0)\geq0$, and consequently 
$(\varDelta^ps)(n)>0$ for every $n\in\Bbb Z_+$. By applying 
Lemma $4.4$ in case $p\in\Bbb Z_+$, since $s=
\varSigma^p(\varDelta^ps)$ we conclude that $s$ is strictly increasing. 
Hence $s(n)>0$ for every $n\in\Bbb Z_+$. From Lemma $4.4$ we 
also derive that if $p\in\Bbb Z_+$, then $s$ is convex.
\newline
We have thus proved $(4.7.1)$ and $(4.7.2)$, plus the final 
assertion.
\comment
$\lim\limits_{n\rightarrow+\infty}\frac{s(n)}{n^p}
=+\infty$
\endcomment
\newline
Now we prove $(4.7.3)$. If $p=0$, the desired result holds as \,
$+\infty=\lim\limits_{n\rightarrow+\infty}(\varDelta^0s)(n)\mathbreak=\lim\limits_
{n\rightarrow+\infty}s(n)$. Thus, let us assume $p\in\Bbb Z_+$. 
Since $\varDelta^ps$ is nondecreasing, from Theorem 
$4.3$ it follows that, for each $n\in\Bbb Z_+$, we have
$$
\multline
s(n)=\bigl(\varSigma^p(\varDelta^ps)\bigr)(n)\geq\frac1{p+1}
\binom{n+p}p(\varDelta^ps)(n)+\frac p{p+1}\binom{n+p}ps(0)
\\\geq\frac1{p+1}\binom{n+p}p(\varDelta^ps)(n)
\endmultline
$$
and consequently
$$
\frac{s(n)}{n^p}\geq\frac{\prod\limits_{k=1}^p(n+k)}
{(p+1)!\,n^p}\cdot(\varDelta^p s)(n)\geq\frac{(\varDelta^ps)(n)}{(p+1)!}.
$$
Since $\lim\limits_{n\rightarrow+\infty}(\varDelta^ps)(n)=+\infty$, we 
obtain the desired result.
\newline
We prove $(4.7.4)$. Since $s(n)\geq0$ for every $n\in\Bbb N$, it 
suffices to prove that the sequence ${\left(\frac{s(n)}{n^{p+1}}\right)}_
{n\in\Bbb Z_+}$ is bounded from above. By Remark $3.4$, 
the sequence ${\left(\frac{(\varDelta^ps)(n)}n\right)}_{n\in\Bbb Z_+}$ is 
bounded from above, which is the desired result if $p=0$. Now suppose 
$p\in\Bbb Z_+$. From Theorem $4.3$ it follows that
$$
\frac{s(n)}{n^{p+1}}=\frac{\bigl(\varSigma^p(\varDelta^ps)\bigr)(n)}
{n^{p+1}}\leq\frac1{p!}
\raise1.3ex\hbox{$\left(\lower1.3ex\hbox{$\displaystyle{\frac{\prod\limits_
{k=1}^p(n+k)}{n^p}}$}\right)$}\!\left(\frac{(\varDelta^ps)(n)}n\right)
\qquad\text{for every }n\in\Bbb Z_+.
$$
Since $\lim\limits_{n\rightarrow+\infty}
\raise.5ex\hbox{$\biggl(\lower.5ex\hbox{$\frac{\prod\limits_{k=1}^p(n+k)}
{n^p}$}\biggr)$}=1$ and ${\left(\frac{(\varDelta^ps)(n)}n\right)}
_{n\in\Bbb Z_+}$ is bounded from above, we obtain the desired result.
\newline
Now we prove $(4.7.5)$. We prove that $\lim\limits_{n\rightarrow
+\infty}\frac{(\varDelta s)(n)}{s(n)}=0$, which is equivalent to proving that 
$\lim\limits_{n\rightarrow+\infty}\frac{s(n+1)}{s(n)}=1$ by Remark $4.6$. 
If $p=0$, then $s$ is concave by hy\-po\-the\-sis. Furthermore, $s$ is 
strictly increasing and $\lim\limits_{n\rightarrow+\infty}s(n)=+\infty$. Since 
$\varDelta s$ converges by Lemma $4.5$, it follows that $\lim\limits_{
n\rightarrow+\infty}\frac{(\varDelta s)(n)}{s(n)}=0$. Now assume $p\in\Bbb
Z_+$ (and consequently $p-1\in\Bbb N$). From Theorem $4.3$ we 
derive that, for each $n\in\Bbb N$, we have
$$
s(n)=\bigl(\varSigma^p(\varDelta^ps)\bigl)(n)\geq\frac1{p+1}\binom
{n+p}p(\varDelta^ps)(n)\qquad\text{(as $(\varDelta^ps)(0)=s(0)\geq0$)}
$$
and
$$
(\varDelta s)(n)=\bigl(\varSigma^{p-1}(\varDelta^ps)\bigl)(n)\leq
\binom{n+p-1}{p-1}(\varDelta^ps)(n).
$$
Since both $s(n)$ and $(\varDelta s)(n)$ are strictly positive for 
every $n\in\Bbb Z_+$ (see $(4.7.1)$ and $(4.7.2)$), from the two 
inequalities above we conclude that
$$
0<\frac{(\varDelta s)(n)}{s(n)}\leq\frac{\binom{n+p-1}{p-1}(\varDelta^ps)(n)}
{\frac1{p+1}\binom{n+p}p(\varDelta^ps)(n)}=\frac{\frac{(n+p-1)!}{(p-1)!\,n!}}
{\frac1{p+1}\cdot\frac{(n+p)!}{p!\,n!}}=\frac{(p+1)!\,(n+p-1)!}{(p-1)!\,(n+p)!}
=\frac{p(p+1)}{n+p}
$$
for every $n\in\Bbb Z_+$. Now $\lim\limits_{n\rightarrow+\infty}\frac
{p(p+1)}{n+p}=0$ \ yields $\lim\limits_{n\rightarrow+\infty}\frac
{(\varDelta s)(n)}{s(n)}=0$, which is the desired result.
\newline
Finally, $(4.7.6)$ is a consequence of Corollary $3.10$ and Lemma 
$4.5$. The proof is now complete.
\qed
\enddemo
\remark{Remark 4.8}
Under the hypotheses of Theorem $4.7$, for each $j=0,\dots,p$ we 
have
\linebreak
$\varDelta^{p-j}(\varDelta^js)=\varDelta^ps$. Since 
$(\varDelta^js)(0)=s(0)\geq0$, we are enabled to apply Theorem $4.7$ 
to the sequence $\varDelta^js$. Hence:
\roster
\item"(4.8.1)"$(\varDelta^js)(n)>0$ for every $n\in\Bbb Z_+$;
\item"(4.8.2)"$\varDelta^js$ is strictly increasing;
\item"(4.8.3)"$\lim\limits_{n\rightarrow+\infty}\frac{(\varDelta^js)(n)}
{n^{p-j}}=+\infty$ 
(and consequently $\lim\limits_{n\rightarrow+\infty}(\varDelta^js)(n)=
+\infty$);
\item"(4.8.4)"the sequence ${\left(\frac{(\varDelta^js)(n)}{n^{p+1-j}}
\right)}_{n\in\Bbb Z_+}$ is bounded;
\item"(4.8.5)"$\lim\limits_{n\rightarrow+\infty}\frac{(\varDelta^js)(n+1)}
{(\varDelta^js)(n)}=1$ 
(or, equivalently, $\lim\limits_{n\rightarrow+\infty}\frac{(\varDelta^{j+1}s)(n)}
{(\varDelta^js)(n)}=0$).
\endroster
If in addition $j<p$, then the sequence $\varDelta^js$ is also convex.
\endremark
\vskip5pt
We conclude this section with an example of an important sequence 
satisfying the hypotheses of Theorem $4.7$, that is, the sequence of 
the Ces\`aro numbers of order $\alpha$ for $\alpha>0$.
\example{Example 4.9}
Fix $\alpha\in(0,+\infty)$, and consider the sequence $A_\alpha:\Bbb
N\rightarrow\Bbb R$ of the Ces\`aro numbers of order $\alpha$. Then 
$A_\alpha(0)=1>0$. Also, from $(2.5)$ it follows that $A_\alpha=
\varSigma A_{\alpha-1}$, or equivalently $\varDelta A_\alpha=
A_{\alpha-1}$. Hence $\varDelta^kA_\alpha=
A_{\alpha-k}$ for each $k\in\Bbb N$. Now we set
$$
p=\max\{k\in\Bbb N:k<\alpha\}.\tag4.9.1
$$
Then
$$
p=\cases[\alpha]&\text{if }\alpha\notin\Bbb Z_+\\\alpha-1&\text{if }
\alpha\in\Bbb Z_+\endcases\qquad\text{and}\qquad0<\alpha-p
\leq1.
$$
Since $\alpha-p>0$, from $(2.6)$ we derive that the sequence 
$\varDelta^pA_\alpha=A_{\alpha-p}$ is unbounded from above. 
Besides, since $-1<\alpha-p-1\leq0$, it follows that the sequence
$\varDelta(\varDelta^pA_\alpha)=A_{\alpha-p-1}$ is nonincreasing 
(see [Z], III, (1-17)), and consequently $\varDelta^pA_\alpha$ is concave. 
Hence $A_\alpha$ satisfies the hypotheses of Theorem $4.7$ if $p$ 
is as in $(4.9.1)$.
\endexample
Notice that a real sequence need not be infinite of order $\alpha$ 
for some $\alpha\in(0,+\infty)$ in order to satisfy the hypotheses of 
Theorem $4.7$: for instance, the sequence ${\bigl(\log(n+1)\bigr)}_
{n\in\Bbb N}$ of nonnegative real numbers, being concave and 
unbounded from above, satisfies the hypotheses 
of Theorem $4.7$ for $p=0$. Nevertheless, it is infinite of order less 
than $\alpha$ for each $\alpha\in(0,+\infty)$.
\heading{\bf 5. An index of unboundedness from above for a real sequence}
\endheading
\definition{Definition 5.1}
For each real sequence $a:\Bbb N\rightarrow\Bbb R$, we set
$$
\Cal H(a)=\inf\bigl\{m\in\Bbb N:\text{ the sequence }{\tsize\bigl(\frac{a(n)}
{n^m}\bigr)}_{n\in\Bbb Z_+}\text{ is bounded from above}\bigr\}.
$$
\enddefinition
\remark{Remark 5.2}
Let $a:\Bbb N\rightarrow\Bbb R$ be a real sequence. 
We observe that $\Cal H(a)\in\Bbb N\cup\{+\infty\}$ and the infimum above 
is attained if and only if $\Cal H(a)<+\infty$. Also, $\Cal H(a)<+\infty$ 
if and only if the sequence ${\bigr(\frac{a(n)}{n^\alpha}\bigr)}_{n\in\Bbb
Z_+}$ is bounded from above for some $\alpha\in[0,+\infty)$, in which 
case $\Cal H(a)$ is the minimum nonnegative integer $m$ for which 
the sequence ${\bigl(\frac{a(n)}{n^m}\bigr)}_{n\in\Bbb Z_+}$ is 
bounded from above. Moreover, if $\Cal H(a)<+\infty$, then clearly 
${\bigl(\frac{a(n)}{n^m}\bigr)}_{n\in\Bbb Z_+}$ is bounded from above 
for every $m\in\Bbb N_{\Cal H(a)}$ (indeed, for every $m\in[\Cal
H(a),+\infty)$). Notice also that $a$ is bounded from above if and 
only if $\Cal H(a)=0$.
\newline
Finally, we remark that if $s:\Bbb N\rightarrow\Bbb R$ is a real sequence 
satisfying the hypotheses of Theorem $4.7$, then $\Cal H(s)=p+1$. Thus, 
if $s:\Bbb N\rightarrow\Bbb R$ is a real sequence such that $s(0)\geq0$, 
there exists at most one $p\in\Bbb N$ for which the hypotheses of Theorem 
$4.7$ are satisfied.
\endremark
\proclaim{Theorem 5.3}
Let $b:\Bbb N\rightarrow\Bbb R$ be a real sequence, which is 
not bounded from above and satisfies $\Cal H(b)<+\infty$. Then 
$\Cal H(b)\in\Bbb Z_+$. Besides, if we set $p=\Cal H(b)-1$, then 
$p\in\Bbb N$ and $b$ has a majorant $s:\Bbb N\rightarrow\Bbb
R$ such that $s(0)\geq0$ and $\varDelta^ps$ is concave and is not 
bounded from above (which implies that $s$ satisfies $(4.7.1)$--$(4.7.6)$, 
besides being convex if $p\in\Bbb Z_+$---equivalently, if $\Cal H(b)
\in\Bbb N_2$), and moreover \,$\limsup\limits_{n\rightarrow+\infty}\frac
{b(n)}{s(n)}\in\bigl[\frac1{p+1},1\bigr]=\bigl[\frac1{\Cal H(b)},1\bigr]$.
\endproclaim
\demo{Proof}
Let us first notice that $\Cal H(b)\in\Bbb Z_+$---and consequently 
$p\in\Bbb N$---by Remark $5.2$, being $b$ unbounded from above. 
By going to the sequence
$$
\tilde b:\Bbb N\ni n\longmapsto\cases b(n)&\text{if }n\in\Bbb Z_+\\
0&\text{if }n=0\endcases\in\Bbb R
$$
if $b(0)<0$, 
it is not restrictive to assume that $b(0)\geq0$. Now let $a:\Bbb N
\rightarrow\Bbb R$ be the sequence defined by
$$
a(n)=\frac{(p+1)b(n)}{\binom{n+p}p}-p\hskip.5pt b(0)\qquad\text
{for every }n\in\Bbb N.
$$
Then $a(0)=b(0)$. Furthermore, since $\lim\limits_{n\rightarrow
+\infty}\frac{\binom{n+p}p}{n^p}=\lim\limits_{n\rightarrow+\infty}
\raise.5ex\hbox{$\biggl(\lower.5ex\hbox{$\frac{\prod\limits_{k=1}^p(n+k)}
{p!n^p}$}\biggr)$}=\frac1{p!}$, and the sequence ${\bigl(\frac{b(n)}
{n^{p+1}}\bigr)}_{n\in\Bbb Z_+}$ is bounded from above, whereas 
${\bigl(\frac{b(n)}{n^p}\bigr)}_{n\in\Bbb Z_+}$ is not (as $\Cal H(b)=
p+1$), it follows that the sequence ${\bigl(\frac{a(n)}n\bigr)}_{n\in\Bbb Z_+}$ 
is bounded from above, and $a$ is not. Hence $a$ has a concave majorant 
by Proposition $3.3$, and consequently (see Remark $3.5$) has a 
least concave majorant. Let $c:\Bbb N\rightarrow\Bbb R$ denote the 
least concave majorant of $a$. Then from $(3.6.1)$ we obtain $c(0)=
a(0)=b(0)$. Moreover, from Theorem $3.9$ we conclude that $c$ is 
strictly increasing, $\lim\limits_{n\rightarrow+\infty}c(n)=+\infty$, and 
$\limsup\limits_{n\rightarrow+\infty}\frac{a(n)}{c(n)}=1$. Hence 
$\lim\limits_{n\rightarrow+\infty}\frac{p\hskip.2pt b(0)}{c(n)}=0$, and 
consequently
$$
\limsup_{n\rightarrow+\infty}\,\frac{b(n)}{\frac1{p+1}\binom{n+p}pc(n)}
=1.\tag5.3.1
$$
Now let $s\in\Bbb R^\Bbb N$ be defined by $s=\varSigma^pc$. We prove 
that $s$ is a majorant of $b$. Since $c$ is concave and nondecreasing, 
from Theorem $4.3$ we derive that, for each $n\in\Bbb N$, we have
$$
\multline
s(n)=(\varSigma^pc)(n)\geq\frac1{p+1}\binom{n+p}pc(n)+\frac p{p+1}
\binom{n+p}pc(0)\\=\frac1{p+1}\binom{n+p}pc(n)+\frac p{p+1}
\binom{n+p}pb(0)\\\geq\frac1{p+1}\binom{n+p}pa(n)+\frac p{p+1}
\binom{n+p}pb(0)\\=\frac1{p+1}\binom{n+p}p\left(\frac{(p+1)b(n)}
{\binom{n+p}p}-p\hskip.5pt b(0)\right)+\frac p{p+1}\binom{n+p}pb(0)\\
=b(n)-\frac p{p+1}\binom{n+p}pb(0)+\frac p{p+1}\binom{n+p}pb(0)=b(n),
\endmultline
$$
which is the desired result. We also remark that $s(0)=c(0)=b(0)\geq0$, 
and $\varDelta^ps=c$ is concave and is not bounded from above.
\newline
Now it remains to prove that $\limsup\limits_{n\rightarrow+\infty}\frac{b(n)}
{s(n)}\in\bigl[\frac1{p+1},1\bigr]$. From Theorem $4.7$ it follows that 
$(\varSigma^pc)(n)=s(n)>0$ for every $n\in\Bbb Z_+$. Then, since $s$ 
is a majorant of $b$, we 
clearly have $\limsup\limits_{n\rightarrow+\infty}\frac{b(n)}{s(n)}\leq1$. 
Finally, we prove that $\limsup\limits_{n\rightarrow+\infty}\frac{b(n)}
{s(n)}\geq\frac1{p+1}$. Since $c$ is strictly increasing, $c(0)\geq0$ yields 
$c(n)>0$ for every $n\in\Bbb Z_+$. Then
$$
\frac{b(n)}{s(n)}=\left(\frac{b(n)}{\frac1{p+1}\binom{n+p}pc(n)}\right)
\left(\frac{\frac1{p+1}\binom{n+p}pc(n)}{(\varSigma^pc)(n)}\right)
\qquad\text{for every }n\in\Bbb Z_+.\tag5.3.2
$$
Since $c$ is concave and nondecreasing, from Theorem $4.3$ we 
conclude that
$$
\frac1{p+1}\leq\frac{\frac1{p+1}\binom{n+p}pc(n)}{(\varSigma^pc)(n)}
\qquad\text{for every }n\in\Bbb Z_+.\tag5.3.3
$$
Now the desired result is a consequence of $(5.3.1)$, $(5.3.2)$ and 
$(5.3.3)$. The proof is thus complete.
\qed
\enddemo
We remark that, by virtue of Theorem $5.3$, any real sequence $b$ 
which is unbounded from above and such that $\Cal H(b)<+\infty$, 
has a majorant $s$ which enjoys the good properties of Theorem 
$4.7$ for $p=\Cal H(b)-1$, satisfies $\Cal H(s)=p+1=\Cal H(b)$ 
(see Remark $5.2$) and is not infinite of higher order than $b$ 
(equivalently, has a subsequence which is infinite of the same 
order as the corresponding subsequence of $b$). Thus, 
in some sense, $s$ is not ''too far'' from $b$.
\proclaim{Proposition 5.4}
Let $a:\Bbb N\rightarrow\Bbb R$ be a real sequence, and $q\in\Bbb
Z_+$ be such that $\varDelta^qa\in\ell_1$. Then $\Cal H(a)\leq q-1$.
\endproclaim
\demo{Proof}
If we set $M={\|\varDelta^qa\|}_{\ell_1}$, for each $n\in\Bbb N$ we have
$$
\multline
(\varDelta^{q-1}a)(n)\leq|(\varDelta^{q-1}a)(n)|=|(\varSigma\varDelta^q
a)(n)|=\left|\,\sum_{k=0}^n(\varDelta^qa)(n)\right|\\\leq\sum_{k=0}^n
|(\varDelta^qa)(n)|\leq M=MA_0(n).
\endmultline
$$
Since the linear operator $\varSigma$ preserves inequalities, and 
consequently $\varSigma^{q-1}$ also does, from $(2.5)$ we conclude 
that
$$
a(n)=(\varSigma^{q-1}\varDelta^{q-1}a)(n)\leq M
(\varSigma^{q-1}A_0)(n)=MA_{q-1}(n)=M\binom{n+q-1}n
$$
for each $n\in\Bbb N$. Since $\lim\limits_{n\rightarrow+\infty}
\frac1{n^{q-1}}\,\binom{n+q-1}n=\frac1{(q-1)!}$\,, and consequently the 
sequence ${\Bigl(\frac1{n^{q-1}}\,\binom{n+q-1}n\Bigr)}_{n\in\Bbb Z_+}$ 
is bounded, it follows that the sequence ${\bigl(\frac{a(n)}{n^{q-1}}\bigr)}_
{n\in\Bbb Z_+}$ is bounded from above. Hence $\Cal H(a)\leq{q-1}$.
\qed
\enddemo
\remark{Remark 5.5}
Let $a:\Bbb N\rightarrow\Bbb R$ be a real sequence. If $\varDelta^qa\in
\ell_1$ for some $q\in\Bbb N$, since $\varDelta(\ell_1)\subseteq\ell_1$ 
it follows that $\varDelta^ka\in\ell_1$ for each $k\in\Bbb N_q$.
\endremark
\remark{Remark 5.6}
If a real sequence $a$ is unbounded from above, and $q\in\Bbb Z_+$ 
is such that $\varDelta^qa\in\ell_1$, then from Proposition $5.4$ it 
follows that $q\geq2$ (as $\Cal H(a)\geq1$).
\endremark
\heading{\bf6. A uniform ergodic theorem for N\"orlund means}
\endheading
We begin with a result relating several properties of the sequence 
of the norms of the iterates of a bounded linear operator.
\proclaim{Theorem 6.1}
Let $X$ be a complex nonzero Banach space, and $T\in L(X)$. Then 
the following conditions are equivalent:
\roster
\item"(6.1.1)"$\Cal H\bigl({({\|T^n\|}_{L(X)}})_{n\in\Bbb N}\bigr)<+\infty$;
\item"(6.1.2)"there exists a sequence $b$ of strictly positive real 
numbers such that $\Cal H(b)<+\infty$, $\lim\limits_{n\rightarrow+\infty}
b(n)=+\infty$, and $\lim\limits_{n\rightarrow+\infty}\frac{{\|T^n\|}_{L(X)}}
{b(n)}=0$;
\item"(6.1.3)"there exists a sequence $s$ of strictly positive real 
numbers such that $\varDelta^ps$ is concave and unbounded from 
above for some $p\in\Bbb N$ (which implies that $s$ satisfies 
$(4.7.2)$--$(4.7.6)$, besides being convex if $p\in\Bbb Z_+$), and 
$\lim\limits_{n\rightarrow+\infty}\frac{{\|T^n\|}_{L(X)}}{s(n)}=0$;
\item"(6.1.4)"there exists a nondecreasing sequence $s$ of strictly 
positive real numbers such that $\lim\limits_{n\rightarrow+\infty}
s(n)=+\infty$, $\lim\limits_{n\rightarrow+\infty}\frac{s(n+1)}{s(n)}
=1$, $\varDelta^qs\in\ell_1$ for some $q\in\Bbb N_2$, and
\linebreak
$\lim\limits_{n\rightarrow+\infty}\frac{{\|T^n\|}_{L(X)}}{s(n)}=0$.
\endroster
Furthermore, if $b:\Bbb N\rightarrow\Bbb R$ is a sequence of 
strictly positive real numbers satisfying $(6.1.2)$, then $\Cal H(b)\in
\Bbb Z_+$, and a sequence $s:\Bbb N\rightarrow\Bbb R$ of strictly 
positive real numbers can be chosen so that $(6.1.3)$ is satisfied 
for $p=\Cal H(b)-1$, $s(n)\geq b(n)$ for each $n\in\Bbb N$, and 
moreover \,$\limsup\limits_{n\rightarrow+\infty}\frac{b(n)}{s(n)}\in
\bigl[\frac1{\Cal H(b)},1\bigr]$.
\newline
Finally, the equivalent conditions $(6.1.1)$--$(6.1.4)$ imply the 
following:
\roster
\item"(6.1.5)"$r(T)\leq1$.
\endroster
\endproclaim
\demo{Proof}
We begin by proving that $(6.1.1)$ implies $(6.1.2)$. If $\Cal
H\bigl({({\|T^n\|}_{L(X)}})_{n\in\Bbb N}\bigr)<+\infty$, it suffices to 
define $b:\Bbb N\rightarrow\nomathbreak\Bbb R$ as 
follows: $b(n)={(n+1)}^{\Cal H\text{\eightpoint(}{({\|T^k\|}_{L(X)})}_
{\!k\in\Bbb N}\text{\,\eightpoint)}+1}$ for every $n\in\Bbb N$. Indeed, 
$b(n)>0$ 
for every $n\in\Bbb N$, $\lim\limits_{n\rightarrow+\infty}\frac
{{\|T^n\|}_{L(X)}}{b(n)}=\lim\limits_{n\rightarrow+\infty}\left(\frac
{{\|T^n\|}_{L(X)}}{{(n+1)}^{\Cal H{\ssize(}{({\|T^k\|}_{L(X)})}_
{\!k\in\Bbb N}\,{\ssize)}}}\right)\bigl(\frac1{n+1}\bigr)=0$, and 
$\Cal H(b)=\Cal H\bigl({({\|T^n\|}_{L(X)}})_{n\in\Bbb N}\bigr)+1<+\infty$.
\newline
Now let us assume that condition $(6.1.2)$ is satisfied by a sequence 
$b$ of strictly positive real numbers. Since $\Cal H(b)<+\infty$ and $b$ 
is unbounded from above, from Theorem $5.3$ it follows that $\Cal
H(b)\in\Bbb Z_+$. Furthermore, if we set $p=\Cal H(b)-1$ (which 
gives $p\in\Bbb N$), then $b$ has a majorant $s$ 
such that $\varDelta^ps$ is concave and is not bounded from above, 
and besides \,$\limsup\limits_{n\rightarrow+\infty}\frac{b(n)}{s(n)}\in
\bigl[\frac1{\Cal H(b)},1\bigr]$. Notice also that 
$s(n)\geq b(n)>0$ for each $n\in\Bbb N$. Finally, since \ $\frac{{\|T^n\|}_
{L(X)}}{s(n)}\leq\frac{{\|T^n\|}_{L(X)}}{b(n)}$ \ for each $n\in\Bbb N$, 
we derive that $\lim\limits_{n\rightarrow+\infty}\frac{{\|T^n\|}_{L(X)}}
{s(n)}=0$. Hence condition $(6.1.3)$ is satisfied by $s$ for $p=\Cal
H(b)-1$.
\newline
From Theorem $4.7$ it follows that $(6.1.3)$ implies $(6.1.4)$. Now we 
prove that $(6.1.4)$ implies $(6.1.1)$. Let $s:\Bbb N\rightarrow\Bbb R$ 
be a nondecreasing sequence of strictly positive real numbers which 
satisfies $(6.1.4)$. Then $\Cal H(s)\leq q-1$ by Proposition $5.4$. 
Also, the sequence ${\left(\frac{{\|T^n\|}_{L(X)}}{s(n)}\right)}_{n\in\Bbb N}$ 
is bounded, and consequently $\Cal H\bigl({({\|T^n\|}_{L(X)}})_
{n\in\Bbb N}\bigr)\leq\Cal H(s)<+\infty$.
\par
We have thus proved equivalence of conditions $(6.1.1)$--$(6.1.4)$, 
as well as the subsequent claim. It remains to prove that if the 
equivalent conditions $(6.1.1)$--$(6.1.4)$ are satisfied, then $r(T)
\leq1$, which follows from Remark $2.8$.
\qed
\enddemo
\comment
We begin with the following corollary of Theorem $5.3$.
\proclaim{Corollary y.1}
Let $X$ be a complex nonzero Banach space, $T\in L(X)$ with 
$r(T)=1$, and $b:\Bbb N\rightarrow\Bbb R$ be a sequence of strictly 
positive real numbers, such that $\lim\limits_{n\rightarrow+\infty}\frac
{{\|T^n\|}_{L(X)}}{b(n)}=0$ and $\Cal H(b)<+\infty$. Then $\lim\limits_
{n\rightarrow+\infty}b(n)=+\infty$. Furthermore, if we set $p=\Cal
H(b)-1$, then $p\in\Bbb N$ and $b$ has a majorant $s:\Bbb
N\rightarrow\Bbb R$ such that $s(0)
>0$ and $\varDelta^ps$ is concave and is not bounded from above 
(which implies that $s$ satisfies $(4.7.1)$--$(4.7.6)$, besides being 
convex if $p\in\Bbb Z_+$), and moreover \,$\limsup\limits_
{n\rightarrow+\infty}\frac{b(n)}{s(n)}\in\bigl[\frac1{p+1},1\bigr]$. Finally, 
$\lim\limits_{n\rightarrow+\infty}\frac{{\|T^n\|}_{L(X)}}{s(n)}=0$.
\endproclaim
\demo{Proof}
Since ${\|T^n\|}_{L(X)}\geq r(T^n)=1$ for every $n\in\Bbb N$, we have
$$
b(n)=\frac{b(n)}{{\|T^n\|}_{L(X)}}\cdot{\|T^n\|}_{L(X)}\geq\frac{b(n)}
{{\|T^n\|}_{L(X)}}@>>n\rightarrow\infty>+\infty.
$$
Hence $\lim\limits_{n\rightarrow+\infty}b(n)=+\infty$. By applying 
Theorem $5.3$, we derive that $p\in\Bbb N$ and $b$ has a 
majorant $s$ such that $s(0)>0$ (as $b(0)>0$), $\varDelta^ps$ is 
concave and is not bounded from above, and \,$\limsup\limits_
{n\rightarrow+\infty}\frac{b(n)}{s(n)}\in\bigl[\frac1{p+1},1\bigr]$. Finally, 
since $s(n)\geq b(n)$ for every $n\in\Bbb N$ and $\lim\limits_
{n\rightarrow+\infty}\frac{{\|T^n\|}_{L(X)}}{b(n)}=0$, it follows that 
$\lim\limits_{n\rightarrow+\infty}\frac{{\|T^n\|}_{L(X)}}{s(n)}=0$, 
which finishes the proof.
\qed
\enddemo
\remark{Remark x.2}
Under the hypotheses of Corollary $y.1$, the sequence $b$ has a 
majorant $s$ which enjoys the good properties of Theorem 
$4.7$ for $p=\Cal H(b)-1$, satisfies $\Cal H(s)=p+1=\Cal H(b)$ 
(see Remark $5.2$) and is not infinite of higher order than $b$. Thus, 
in some sense, $s$ is not ''too far'' from $b$.
\endremark
\remark{Remark z.3}
If $T$ is a bounded linear operator on a complex 
nonzero Banach space $X$, such that $r(T)=1$, then a real sequence 
$b$ satisfying the hypotheses of Corollary $y.1$ exists if and only 
if $\Cal H\bigl({({\|T^n\|}_{L(X)}})_{n\in\Bbb N}\bigr)<+\infty$. Indeed, 
if $b:\Bbb N\rightarrow\Bbb R$ is as in Corollary $y.1$, then the 
sequence ${\left(\frac{{\|T^n\|}_{L(X)}}{b(n)}\right)}_{n\in\Bbb N}$ 
is bounded, and consequently $\Cal H\bigl({({\|T^n\|}_{L(X)}})_
{n\in\Bbb N}\bigr)\leq\Cal H(b)<+\infty$. Conversely, if $\Cal
H\bigl({({\|T^n\|}_{L(X)}})_{n\in\Bbb N}\bigr)<+\infty$, it suffices to 
define $b:\Bbb N\rightarrow\nomathbreak\Bbb R$ as follows: $b(n)=
{(n+1)}^{\Cal H\text{\eightpoint(}{({\|T^k\|}_{L(X)})}_{\!k\in\Bbb N}
\text{\,\eightpoint)}+1}$ for every $n\in\Bbb N$. Then $b(n)>0$ 
for every $n\in\Bbb N$, $\lim\limits_{n\rightarrow+\infty}\frac
{{\|T^n\|}_{L(X)}}{b(n)}=\lim\limits_{n\rightarrow+\infty}\left(\frac
{{\|T^n\|}_{L(X)}}{{(n+1)}^{\Cal H{\ssize(}{({\|T^k\|}_{L(X)})}_
{\!k\in\Bbb N}\,{\ssize)}}}\right)\bigl(\frac1{n+1}\bigr)=0$, and 
$\Cal H(b)=\Cal H\bigl({({\|T^n\|}_{L(X)}})_{n\in\Bbb N}\bigr)+1<+\infty$.
\endremark
\endcomment
\remark{Remark 6.2}
If $T$ is a bounded linear operator on a complex 
nonzero Banach space $X$, such that $r(T)<1$, then $\lim\limits_
{n\rightarrow+\infty}{\|T^n\|}_{L(X)}=0$. Consequently, the 
sequence ${{({\|T^n\|}_{L(X)}})}_{n\in\Bbb N}$ is bounded, which 
gives $\Cal H\bigl({({\|T^n\|}_{L(X)}})_{n\in\Bbb N}\bigr)=0$.
\newline
However, condition $(6.1.5)$ is not equivalent to 
$(6.1.1)$--$(6.1.4)$. Indeed, the following example 
shows that a bounded linear operator $T$ on a complex nonzero 
Banach space $X$, such that $r(T)=1$, need not satisfy $\Cal
H\bigl({({\|T^n\|}_{L(X)}})_{n\in\Bbb N}\bigr)<+\infty$.
\endremark
\example{Example 6.3}
Let us consider the complex Hilbert space $\ell_2$ and the 
unilateral weighted shift operator $T:\ell_2\rightarrow\ell_2$ defined by
$$
Tx=\sum_{n=0}^{+\infty}e^{\frac1{\sqrt{(n+1)}}}\,x(n)\,e_{n+1}
\qquad\text{for every }x\in\ell_2,
$$
where $\{e_n:n\in\Bbb N\}$ denotes the canonical orthonormal basis 
of $\ell_2$. Then $T\in L(\ell_2)$. Besides (see [Ha], Solution $77$), 
for each $k\in\Bbb Z_+$ we have
$$
{\|T^k\|}_{L(\ell_2)}=\sup\Biggl\{\ \prod_{j=1}^ke^{\frac1{\sqrt{n+j}}}:
n\in\Bbb N\Biggr\}=\sup\Biggl\{e^{\,\sum\limits_{j=1}^k\!\frac1
{\sqrt{n+j}}}:n\in\Bbb N\Biggr\}=e^{\,\sum\limits_{j=1}^k\!\frac1
{\sqrt j}}.
$$
Since $\frac1{\sqrt j}\rightarrow0$ as $j\rightarrow+\infty$, from 
the classical Ces\`aro means theorem we conclude that
$$
\root k\of{{{\|T^k\|}_{L(\ell_2)}}}=e^{\frac1k\!\sum\limits_{j=1}^k\!\frac1
{\sqrt j}}@>>k\rightarrow+\infty>e^0=1,
$$
and consequently $r(T)=1$.
\newline
Now fix $\alpha\in(0,+\infty)$. Since for each $j\in\Bbb Z_+$ we have 
\,$\frac1{\sqrt j}\geq\frac1{\sqrt x}$ \,for every $x\in[j,j+1]$, and 
consequently \,$\frac1{\sqrt j}\geq\int_j^{j+1}\!\frac1{\sqrt x}\,dx$, 
it follows that
$$
\sum_{j=1}^k\frac1{\sqrt j}\geq\sum_{j=1}^k\,\int_j^{j+1}\!\frac1{\sqrt x}\,dx
=\int_1^{k+1}\frac1{\sqrt x}\ dx=2\sqrt{k+1}-2\qquad\text{ for every }k\in
\Bbb Z_+.
$$
Then
$$
\frac{{\|T^k\|}_{L(\ell_2)}}{k^\alpha}=e^{\,\sum\limits_{j=1}^k\!\frac1
{\sqrt j}-\alpha\log k}\geq e^{2\sqrt{k+1}-\alpha\log k-2}
@>>k\rightarrow+\infty>+\infty.
$$
Hence the sequence ${\Bigl(\frac{{\|T^n\|}_{L(X)}}{n^\alpha}\Bigr)}_
{n\in\Bbb Z_+}$ is bounded from above for no $\alpha\in(0,+\infty)$, 
that is, $\Cal H\bigl({({\|T^n\|}_{L(X)}})_{n\in\Bbb N}\bigr)=+\infty$.
\endexample
\proclaim{Lemma 6.4}
Let $\Cal A$ be an algebra with identity $\bold1_\Cal A$ over $\Bbb K$, 
$\tau\in\Cal A$, $a\in\Bbb K^\Bbb N$. Then for each $m\in\Bbb Z_+$ 
and each $n\in\Bbb N$ we have
$$
\multline
\left(\,\sum_{k=0}^na(n-k)\,\tau^k\right){(\bold1_\Cal A-\tau)}^m\\
={(-1)}^m\sum_{k=0}^{n+m}(\varDelta^ma)(n+m-k)\,\tau^k+
\sum_{j=0}^{m-1}{(-1)}^j(\varDelta^ja)(n+j+1){(\bold1_\Cal A-\tau)}^{m-1-j}.
\endmultline
$$
\endproclaim
\demo{Proof}
We proceed by induction on $m$. We set
$$
\multline
S=\biggl\{m\in\Bbb Z_+:\biggl(\ \sum_{k=0}^na(n-k)\,\tau^k\biggr)
{(\bold1_\Cal A-\tau)}^m={(-1)}^m\sum_{k=0}^{n+m}(\varDelta^ma)(n+m-k)
\,\tau^k\\+\sum_{j=0}^{m-1}{(-1)}^j(\varDelta^ja)(n+j+1){(\bold1_\Cal A-\tau)}^
{m-1-j}\text{ \,for every }n\in\Bbb N\biggr\}.
\endmultline
$$
Since for each $n\in\Bbb N$ we have
$$
\multline
\left(\,\sum_{k=0}^na(n-k)\,\tau^k\right)(\bold1_\Cal A-\tau)=\sum_{k=0}^n
a(n-k)\,\tau^k-\sum_{k=0}^na(n-k)\,\tau^{k+1}\\=\sum_{k=0}^na(n-k)\,\tau^k
-\sum_{k=1}^{n+1}a(n+1-k)\,\tau^k\\=-\sum_{k=0}^na(n+1-k)\,\tau^k+
\sum_{k=0}^na(n-k)\,\tau^k-a(0)\,\tau^{n+1}+a(n+1)\,\bold1_\Cal A\\
=-\left(\,\sum_{k=0}^n(\varDelta a)(n+1-k)\,\tau^k+(\varDelta a)(0)\,\tau^
{n+1}\right)+(\varDelta^0a)(n+1)\,\bold1_\Cal A\\=-\left(\,\sum_{k=0}^
{n+1}(\varDelta a)(n+1-k)\,\tau^k\right)+(\varDelta^0a)(n+1)
{(\bold1_\Cal A-\tau)}^0,
\endmultline
$$
it follows that $1\in S$.
\newline
Now let $m\in S$. Then, since $(\varDelta^ma)(0)=a(0)=(\varDelta^{m+1}a)(0)$, 
for each $n\in\Bbb N$ we have
$$
\multline
\left(\,\sum_{k=0}^na(n-k)\,\tau^k\right){(\bold1_\Cal A-\tau)}^{m+1}=
\left(\biggl(\ \sum_{k=0}^na(n-k)\,\tau^k\biggr){(\bold1_\Cal A-\tau)}^m\right)
(\bold1_\Cal A-\tau)\\\allowdisplaybreak=\Biggl({(-1)}^m\sum_{k=0}^{n+m}
(\varDelta^ma)(n+m-k)\,\tau^k+\sum_{j=0}^{m-1}{(-1)}^j(\varDelta^ja)(n+j+1)
{(\bold1_\Cal A-\tau)}^{m-1-j}\Biggr)(\bold1_\Cal A-\tau)\\\allowdisplaybreak
={(-1)}^m\sum_{k=0}^{n+m}(\varDelta^ma)(n+m-k)\,\tau^k+{(-1)}^{m+1}\sum_
{k=0}^{n+m}(\varDelta^ma)(n+m-k)\,\tau^{k+1}\\\allowdisplaybreak+\sum_
{j=0}^{m-1}{(-1)}^j(\varDelta^ja)(n+j+1){(\bold1_\Cal A-\tau)}^{m-j}\\
\allowdisplaybreak={(-1)}^{m+1}\left(\,\sum_{k=1}^{n+m+1}(\varDelta^ma)
(n+m+1-k)\,\tau^k-\sum_{k=0}^{n+m}(\varDelta^ma)(n+m-k)\,\tau^k\right)\\
\allowdisplaybreak+\sum_{j=0}^{m-1}{(-1)}^j(\varDelta^ja)(n+j+1)
{(\bold1_\Cal A-\tau)}^{m-j}\\\allowdisplaybreak={(-1)}^{m+1}
\biggr(\,\sum_{k=0}^{n+m}(\varDelta^{m+1}a)(n+m+1-k)\,\tau^k+(\varDelta^ma)(0)
\,\tau^{n+m+1}-(\varDelta^ma)(n+m+1)\,\bold1_\Cal A\biggl)\\\allowdisplaybreak
+\sum_{j=0}^{m-1}{(-1)}^j(\varDelta^ja)(n+j+1){(\bold1_\Cal A-\tau)}^{m-j}\\
\allowdisplaybreak={(-1)}^{m+1}\biggr(\,\sum_{k=0}^{n+m}(\varDelta^{m+1}a)
(n+m+1-k)\,\tau^k+(\varDelta^{m+1}a)(0)\,\tau^{n+m+1}\biggl)\\\allowdisplaybreak
+{(-1)}^m(\varDelta^ma)(n+m+1)\,\bold1_\Cal A+\sum_{j=0}^{m-1}{(-1)}^j
(\varDelta^ja)(n+j+1){(\bold1_\Cal A-\tau)}^{m-j}\\\allowdisplaybreak=
{(-1)}^{m+1}\biggr(\,\sum_{k=0}^{n+m+1}(\varDelta^{m+1}a)(n+m+1-k)
\,\tau^k\biggl)+\sum_{j=0}^m{(-1)}^j
(\varDelta^ja)(n+j+1){(\bold1_\Cal A-\tau)}^{m-j},
\endmultline
$$
from which we conclude that $m+1\in S$. The proof is now complete.
\qed
\enddemo
\proclaim{Lemma 6.5}
Let $s\in\Bbb K^\Bbb N$ be an eventually nonzero sequence, such that 
$\lim\limits_{n\rightarrow+\infty}\frac{s(n+1)}{s(n)}=1$. Then for each 
$k\in\Bbb Z_+$ we have $\lim\limits_{n\rightarrow+\infty}\frac
{(\varDelta^ks)(n)}{s(n)}=0$ and $\lim\limits_{n\rightarrow+\infty}\frac
{s(n+k)}{s(n)}=1$.
\endproclaim
\demo{Proof}
We begin by proving that $\lim\limits_{n\rightarrow+\infty}\frac
{(\varDelta^ks)(n)}{s(n)}=0$ for each $k\in\Bbb Z_+$, proceeding by induction. 
From Remark $4.6$ it follows that $\lim\limits_{n\rightarrow+\infty}\frac
{(\varDelta s)(n)}{s(n)}=0$. Now let $k\in\Bbb Z_+$ be such that $\lim\limits_
{n\rightarrow+\infty}\frac{(\varDelta^ks)(n)}{s(n)}=0$. Since for each $n\in\Bbb
N$ such that $s(n)\neq0$---and therefore for sufficiently large $n$---we 
have
$$
\frac{(\varDelta^{k+1}s)(n)}{s(n)}=\frac{(\varDelta^ks)(n)}{s(n)}-\frac
{(\varDelta^ks)(n-1)}{s(n-1)}\cdot\frac{s(n-1)}{s(n)},
$$
we conclude that $\lim\limits_{n\rightarrow+\infty}\frac{(\varDelta^{k+1}s)(n)}
{s(n)}=0$, which gives the desired result.
\newline
Now, in order to finish the proof of the lemma, it suffices to observe that for 
each $k\in\Bbb Z_+$ we have
$$
\frac{s(n+k)}{s(n)}=\prod_{j=0}^{k-1}\frac{s(n+j+1)}
{s(n+j)}@>>n\rightarrow+\infty>1.
$$
\qed
\enddemo
\definition{Definition 6.6}
If $X$ is a normed space and $T\in L(X)$, let $\Cal M_T:\Bbb N\rightarrow
\Bbb R$ be the real sequence defined by
$$
\Cal M_T(n)=\max\bigl\{{\|T^k\|}_{L(X)}:k=0,\dots,n\bigr\}\qquad\text{for every }
n\in\Bbb N.
$$
\enddefinition
\proclaim{Theorem 6.7}
Let $X$ be a complex nonzero Banach space, $T\in L(X)$, and $b:\Bbb
N\rightarrow\Bbb R$ be a sequence of strictly positive real numbers, such 
that $\Cal H(b)<+\infty$, $\lim\limits_{n\rightarrow+\infty}b(n)=+\infty$ and 
$\lim\limits_{n\rightarrow+\infty}\frac{{\|T^n\|}_{L(X)}}{b(n)}=0$. Then $r(T)
\leq1$. Furthermore, if $s:\Bbb N\rightarrow\Bbb R$ is any nondecreasing 
sequence of strictly positive real numbers, such that $\lim\limits_
{n\rightarrow+\infty}s(n)=+\infty$, $\lim\limits_{n\rightarrow+\infty}
\frac{s(n+1)}{s(n)}=1$, $\varDelta^qs\in\ell_1$ for some 
$q\in\Bbb N_2$, and the sequence ${\bigl(\frac{b(n)}{s(n)}\bigr)}_
{\!n\in\Bbb N}$ is bounded\footnote{Notice that, by virtue of Theorem 
$6.1$, such a sequence $s$ exists, and can be chosen so that it is not 
infinite of higher order than $b$.}, then $\lim\limits_{n\rightarrow+\infty}\frac
{{\|T^n\|}_{L(X)}}{s(n)}=0$, and the following conditions are equivalent:
\roster
\item"(6.7.1)"the sequence ${\Biggl(\frac{\sum\limits_{k=0}^n(\varDelta s)(n-k)
\,T^k}{s(n)}\Biggr)}_{\!n\in\Bbb N}$ converges in $L(X)$;
\item"(6.7.2)"$1$ is either in $\rho(T)$, or a simple pole of $\goth R_T$;
\item"(6.7.3)"$X=\Cal N(I_X-T)\oplus\Cal R(I_X-T)$;
\item"(6.7.4)"$\Cal R(I_X-T)$ is closed in $X$ and $X=\Cal N(I_X-T)\oplus\Cal
R(I_X-T)$.
\endroster
Finally, if the equivalent conditions $(6.7.1)$--$(6.7.4)$ are 
satisfied and $P\in L(X)$ is such that \ \,$\frac{\sum\limits_{k=0}^n(\varDelta s)
(n-k)\,T^k}{s(n)}\longrightarrow P$ in $L(X)$ as $n\rightarrow+\infty$, then 
$P$ is the projection of $X$ onto $\Cal N(I_X-T)$ along $\Cal R(I_X-T)$.
\endproclaim
\demo{Proof}
We begin by remarking that $r(T)\leq1$ by Theorem $6.1$. Now let 
$s:\Bbb N\rightarrow\Bbb R$ be a nondecreasing 
sequence of strictly positive real numbers, such that $\lim\limits_
{n\rightarrow+\infty}s(n)=+\infty$, $\lim\limits_{n\rightarrow+\infty}
\frac{s(n+1)}{s(n)}=1$, $\varDelta^qs\in\ell_1$ for some 
$q\in\Bbb N_2$, and the sequence ${\bigl(\frac{b(n)}{s(n)}\bigr)}_
{\!n\in\Bbb N}$ is bounded. Clearly, $\lim\limits_{n\rightarrow+\infty}\frac
{{\|T^n\|}_{L(X)}}{b(n)}=0$ yields $\lim\limits_{n\rightarrow+\infty}\frac
{{\|T^n\|}_{L(X)}}{s(n)}=0$. We prove that conditions $(6.7.1)$--$(6.7.4)$ 
are equivalent.
\newline
We first observe that conditions $(6.7.2)$--$(6.7.4)$ are equivalent 
by Theorem $2.1$. Now suppose that the equivalent conditions 
$(6.7.2)$--$(6.7.4)$ are satisfied, and let $P$ denote the projection of 
$X$ onto $\Cal N(I_X-T)$ along $\Cal R(I_X-T)$. Then $P\in
L(X)$. We prove that \ $\frac{\sum\limits_{k=0}^n(\varDelta s)
(n-k)\,T^k}{s(n)}\longrightarrow P$ in $L(X)$ as 
$n\rightarrow+\infty$.
\newline
Since $Tx=x$ for every $x\in\Cal N(I_X-T)$, it follows that $TP=P$, and 
consequently $T^kP=P$ for every $k\in\Bbb N$. Then for each $n\in
\Bbb N$ we have
$$
\multline
\left(\frac{\sum\limits_{k=0}^n(\varDelta s)(n-k)\,T^k}{s(n)}\right)P=
\frac{\sum\limits_{k=0}^n(\varDelta s)(n-k)\,P}{s(n)}=
\left(\frac{\sum\limits_{k=0}^n(\varDelta s)(n-k)}{s(n)}\right)P\\
\allowdisplaybreak=\left(\frac{\sum\limits_{j=0}^n(\varDelta s)(j)}
{s(n)}\right)P=\biggl(\frac{(\varSigma\varDelta s)(n)}{s(n)}\biggr)P
=\biggl(\frac{s(n)}{s(n)}\biggr)P=P.
\endmultline
\tag6.7.5
$$
\comment
Hence $\lim\limits_{n\rightarrow+\infty}\Biggl(\frac{\sum\limits_{k=0}^n
(\varDelta s)(n-k)\,T^k}{s(n)}\Biggr)P=P$ in $L(X)$.
\newline
\endcomment
Now we prove that 
$\Biggl(\frac{\sum\limits_{k=0}^n(\varDelta s)(n-k)\,T^k}{s(n)}\Biggr)(I_X-P)
\longrightarrow0_{L(X)}$ in $L(X)$ as $n\rightarrow+\infty$. Since $T$ 
satisfies the equivalent conditions $(6.7.2)$--$(6.7.4)$, from Theorem $2.1$ 
it follows that $\Cal N\bigl({(I_X-T)}^n\bigr)=\Cal N(I_X-T)$ and $\Cal
R\bigl({(I_X-T)}^n\bigr)=\Cal R(I_X-T)$ for every $n\in\Bbb Z_+$. Then the 
bounded linear operator
$$
A:\Cal R(I_X-T)\ni x\longmapsto{(I_X-T)}^{q-1}x\in\Cal R(I_X-T)
$$
is bijective: indeed, since $q\in\Bbb N_2$ (and so $q-1\in\Bbb
Z_+$), we have
$$
\Cal N(A)=\Cal N\bigl({(I_X-T)}^{q-1}\bigr)\cap\Cal R(I_X-T)=\Cal
N(I_X-T)\cap\Cal R(I_X-T)=\{0_X\},
$$
and
$$
\Cal R(A)=\Cal R\bigl({(I_X-T)}^q\bigr)=\Cal R(I_X-T).
$$
Since $\Cal R(I_X-T)$ is a closed subspace of $X$, and consequently a 
Banach space, it follows that the linear map $A^{-1}:\Cal R(I_X-T)
\longrightarrow\Cal R(I_X-T)$ is bounded. Since $I_X-P$ is the projection 
of $X$ onto $\Cal R(I_X-T)$ along $\Cal N(I_X-T)$, and consequently 
$\Cal R(I_X-P)=\Cal R(I_X-T)$, that is the domain of $A^{-1}$, we can 
define the linear operator
$$
B:X\ni x\longmapsto A^{-1}(I_X-P)x\in X.
$$
We remark that $B\in L(X)$ and
$$
{(I_X-T)}^{q-1}B=I_X-P.\tag6.7.6
$$
By virtue of Lemma $6.4$, for each $n\in\Bbb N$ we have
$$
\multline
\left(\frac{\sum\limits_{k=0}^n(\varDelta s)(n-k)\,T^k}{s(n)}\right)
{(I_X-T)}^{q-1}=\\\frac{{(-1)}^{q-1}\sum\limits_{k=0}^{n+q-1}
(\varDelta^qs)(n+q-1-k)\,T^k+\sum\limits_{j=0}^{q-2}{(-1)}^j
(\varDelta^{j+1}s)(n+j+1){(I_X-T)}^{q-2-j}}{s(n)},
\endmultline
$$
from which we conclude that, for each $n\in\Bbb N$,
$$
\multline
\left\|\left(\frac{\sum\limits_{k=0}^n(\varDelta s)(n-k)\,T^k}{s(n)}\right)
{(I_X-T)}^{q-1}\right\|\\\allowdisplaybreak
\leq\sum_{j=0}^{q-2}\,\frac{|(\varDelta^{j+1}s)(n+j+1)|}{s(n)}\,
{\|I_X-T\|}^{q-2-j}_{L(X)}+\frac{\sum\limits_{k=0}^{n+q-1}|(\varDelta^qs)
(n+q-1-k)|{\|T^k\|}_{L(X)}}{s(n)}\\\allowdisplaybreak
\leq\sum_{j=0}^{q-2}\,\frac{|(\varDelta^{j+1}s)(n+j+1)|}{s(n)}\,{\|I_X-T\|}^{q-2-j}_{L(X)}
+\frac{\Cal M_T(n+q-1)}{s(n)}\sum_{k=0}^{n+q-1}|(\varDelta^qs)(n+q-1-k)|
\\\allowdisplaybreak
=\sum_{j=0}^{q-2}\,\frac{|(\varDelta^{j+1}s)(n+j+1)|}{s(n)}\,{\|I_X-T\|}^{q-2-j}_{L(X)}
+\frac{\Cal M_T(n+q-1)}{s(n)}\sum_{h=0}^{n+q-1}|(\varDelta^qs)(h)|
\\\allowdisplaybreak
\leq\sum_{j=0}^{q-2}\,\frac{|(\varDelta^{j+1}s)(n+j+1)|}{s(n)}\,{\|I_X-T\|}^{q-2-j}_{L(X)}
+\frac{\Cal M_T(n+q-1)}{s(n)}\,{\|\varDelta^qs\|}_{\ell_1}.
\endmultline
\tag6.7.7
$$
For each $j\in\{0,\dots,q-2\}$, we have
$$
\frac{(\varDelta^{j+1}s)(n+j+1)}{s(n)}=\frac{(\varDelta^{j+1}s)(n+j+1)}
{s(n+j+1)}\cdot\frac{s(n+j+1)}{s(n)}\quad
\text{for each }n\in\Bbb N.\tag6.7.8
$$
From Lemma $6.5$ it follows that
$$
\frac{(\varDelta^{j+1}s)(n+j+1)}{s(n+j+1)}@>>n\rightarrow+\infty>0
\qquad\text{and}\qquad\frac{s(n+j+1)}{s(n)}@>>n\rightarrow+\infty>1.\tag6.7.9
$$
Now from $(6.7.8)$ and $(6.7.9)$ we obtain
$$
\lim_{n\rightarrow+\infty}\frac{(\varDelta^{j+1}s)(n+j+1)}{s(n)}=0
\qquad\text{for all }j=0,\dots,q-2,
$$
from which we derive that
$$
\sum_{j=0}^{q-2}\,\frac{|(\varDelta^{j+1}s)(n+j+1)|}{s(n)}\,{\|I_X-T\|}^{q-2-j}_{L(X)}
@>>n\rightarrow+\infty>0.\tag6.7.10
$$
By hypothesis, $s$ is nondecreasing and $s(n)>0$ for each $n\in\Bbb N$. 
Then the sequence ${\bigl(\frac1{s(n)}\bigr)}_{n\in\Bbb N}$ is nonincreasing. 
Since $\lim\limits_{n\rightarrow+\infty}\frac{{\|T^n\|}_{L(X)}}{s(n)}=0=\lim
\limits_{n\rightarrow+\infty}\frac1{s(n)}$ (as $\lim\limits_{n\rightarrow+\infty}
s(n)=+\infty$), from [B], $2.3$ we conclude that $\lim\limits_{n\rightarrow+\infty}
\frac{\Cal M_T(n)}{s(n)}=0$. Since $\lim\limits_{n\rightarrow+\infty}\frac
{s(n+q-1)}{s(n)}=1$ by Lemma $6.5$, we derive that
$$
\frac{\Cal M_T(n+q-1)}{s(n)}=\frac{\Cal M_T(n+q-1)}{s(n+q-1)}\cdot
\frac{s(n+q-1)}{s(n)}@>>n\rightarrow+\infty>0.
$$
This, together with $(6.7.10)$ and $(6.7.7)$, gives
$$
\left(\frac{\sum\limits_{k=0}^n(\varDelta s)(n-k)\,T^k}{s(n)}\right)
{(I_X-T)}^{q-1}@>>n\rightarrow+\infty>0_{L(X)}\qquad\text{in }L(X).
$$
Consequently, by $(6.7.6)$,
$$
\multline
\left(\frac{\sum\limits_{k=0}^n(\varDelta s)(n-k)\,T^k}{s(n)}\right)(I_X-P)
\\=\left(\frac{\sum\limits_{k=0}^n(\varDelta s)(n-k)\,T^k}{s(n)}\right)
{(I_X-T)}^{q-1}B@>>n\rightarrow+\infty>0_{L(X)}\qquad\text{in }L(X).
\endmultline
\tag6.7.11
$$
Now from $(6.7.5)$ and $(6.7.11)$ we conclude that
$$
\multline
\frac{\sum\limits_{k=0}^n(\varDelta s)(n-k)\,T^k}{s(n)}=
\left(\frac{\sum\limits_{k=0}^n(\varDelta s)(n-k)\,T^k}{s(n)}\right)P+
\left(\frac{\sum\limits_{k=0}^n(\varDelta s)(n-k)\,T^k}{s(n)}\right)(I_X-P)
\\=P+\left(\frac{\sum\limits_{k=0}^n(\varDelta s)(n-k)\,T^k}{s(n)}\right)
(I_X-P)@>>n\rightarrow+\infty>P\qquad\text{in }L(X).
\endmultline
$$
We have thus proved that if the equivalent conditions $(6.7.2)$--$(6.7.4)$ 
are satisfied, then the sequence 
${\Biggl(\frac{\sum\limits_{k=0}^n(\varDelta s)(n-k)\,T^k}{s(n)}\Biggr)}_
{\!\!n\in\Bbb N}$ converges in $L(X)$ to the projection of $X$ onto $\Cal
N(I_X-T)$ along $\Cal R(I_X-T)$. It remains to prove that if the 
sequence \ ${\Biggl(\frac{\sum\limits_{k=0}^n(\varDelta s)(n-k)\,T^k}
{s(n)}\Biggr)}_{\!\!n\in\Bbb N}$ \ converges in $L(X)$, then the equivalent 
conditions $(6.7.2)$--$(6.7.4)$ hold.
\newline
Let $P\in L(X)$ be such that \ $\frac{\sum\limits_{k=0}^n(\varDelta s)(n-k)\,T^k}
{s(n)}\longrightarrow P$ in $L(X)$ as $n\rightarrow+\infty$. Since
$$
\frac{\sum\limits_{k=0}^{n+1}(\varDelta s)(n+1-k)\,T^k}{s(n+1)}@>>n\rightarrow
+\infty>P\quad\text{in }L(X)\qquad\text{and}\qquad\frac{s(n+1)}{s(n)}@>>n
\rightarrow+\infty>1,
$$
it follows that
$$
\frac{\sum\limits_{k=0}^{n+1}(\varDelta s)(n+1-k)\,T^k}{s(n)}@>>n\rightarrow
+\infty>P\quad\text{in }L(X),
$$
which in turn yields the following limit in $L(X)$.
$$
\multline
0_{L(X)}=\lim_{n\rightarrow+\infty}\left(\frac{\sum\limits_{k=0}^n(\varDelta s)(n-k)
\,T^k}{s(n)}-\frac{\sum\limits_{k=0}^{n+1}(\varDelta s)(n+1-k)\,T^k}{s(n)}\right)
\\\allowdisplaybreak
=\lim_{n\rightarrow+\infty}\left(\frac{\sum\limits_{k=0}^n(\varDelta s)(n-k)
\,T^k}{s(n)}-\frac{(\varDelta s)(n+1)\,I_X+\sum\limits_{k=1}^{n+1}(\varDelta s)
(n+1-k)\,T^k}{s(n)}\right)
\\\allowdisplaybreak
=\lim_{n\rightarrow+\infty}\left(\frac{\sum\limits_{k=0}^n(\varDelta s)(n-k)
\,T^k}{s(n)}-\frac{(\varDelta s)(n+1)\,I_X+\sum\limits_{k=0}^n(\varDelta s)
(n-k)\,T^{k+1}}{s(n)}\right)\\\allowdisplaybreak=\lim_{n\rightarrow+\infty}
\left(\frac{(I_X-T)\sum\limits_{k=0}^n(\varDelta s)(n-k)\,T^k}{s(n)}-\frac
{(\varDelta s)(n+1)}{s(n)}\,I_X\right)\\\allowdisplaybreak
=\lim_{n\rightarrow+\infty}\,(I_X-T)\!\left(\frac{\sum\limits_{k=0}^n(\varDelta s)(n-k)
\,T^k}{s(n)}\right)
\endmultline
\tag6.7.12
$$
(as $\lim\limits_{n\rightarrow+\infty}\frac{(\varDelta s)(n+1)}{s(n+1)}=0$ by 
Remark $4.6$, being $\lim\limits_{n\rightarrow+\infty}\frac{s(n+1)}{s(n)}=1$, 
and consequently $\frac{(\varDelta s)(n+1)}{s(n)}=\frac{(\varDelta s)(n+1)}{s(n+1)}
\cdot\frac{s(n+1)}{s(n)}\longrightarrow0$ as $n\rightarrow+\infty$).
\vskip5pt
\noindent
Now, for each $n\in\Bbb N$, let $f_n:\Bbb C\longrightarrow\Bbb C$ be the 
polynomial defined by
$$
f_n(z)=\frac{\sum\limits_{k=0}^n(\varDelta s)(n-k)\,z^k}{s(n)}\qquad\text
{for each }z\in\Bbb C.
$$
Since $f_n(1)=\frac{(\varSigma\varDelta s)(n)}{s(n)}=\frac{s(n)}{s(n)}=1$ and 
$f_n(T)=\frac{\sum\limits_{k=0}^n(\varDelta s)(n-k)\,T^k}{s(n)}$ \,for 
every $n\in\Bbb N$, $(6.7.12)$ enables us to apply Theorem $2.2$ 
(together with Remark $2.3$) to the sequence ${(f_n)}_{n\in\Bbb N}$, 
and to conclude that the equivalent conditions $(6.7.2)$--$(6.7.4)$ 
are satisfied. This finishes the proof.
\qed
\enddemo
\remark{Remark 6.8}
Let $X$ be a complex nonzero Banach space, $T\in L(X)$, and $s:\Bbb
N\rightarrow\Bbb R$ be a nondecreasing sequence of strictly positive 
real numbers, such that $\lim\limits_{n\rightarrow+\infty}s(n)=+\infty$, 
$\lim\limits_{n\rightarrow+\infty}\frac{s(n+1)}{s(n)}=1$, $\varDelta^qs\in
\ell_1$ for some $q\in\Bbb N_2$ (which implies $\Cal H(s)\leq q-1$ by 
Proposition $5.4$), and $\lim\limits_{n\rightarrow+\infty}\frac{{\|T^n\|}_
{L(X)}}{s(n)}=0$. Then $r(T)\leq1$ by Theorem $6.1$. Furthermore, from 
Theorem $6.7$ and Theorem $2.4$ we derive that, given any $E\in L(X)$, 
the following two conditions are equivalent:
\roster
\item"(6.8.1)"${\lim\limits_{n\rightarrow+\infty}\left\|\frac{\sum\limits_{k=0}^n
(\varDelta s)(n-k)\,T^k}{s(n)}-E\right\|}_{L(X)}=0$;
\item"(6.8.2)"$\lim\limits_{\lambda\rightarrow1^+}{\|(\lambda-1)\,\goth
R_T(\lambda)-E\|}_{L(X)}=0$.
\endroster
Also, if the equivalent conditions $(6.8.1)$ and $(6.8.2)$ are satisfied, 
then $1$ is either in $\rho(T)$, or a simple pole of $\goth R_T$ (so that 
$\Cal R(I_X-T)$ is closed in $X$ and $X=\Cal N(I_X-T)\oplus\Cal
R(I_X-T)$), and $E$ is the projection of $X$ onto $\Cal N(I_X-T)$ 
along $\Cal R(I_X-T)$.
\endremark
\vskip5pt
The following is a consequence of Theorem $6.7$ and Theorem $4.7$.
\proclaim{Corollary 6.9}
Let $X$ be a complex nonzero Banach space, $T\in L(X)$, and $b:\Bbb
N\rightarrow\Bbb R$ be a sequence of strictly positive real numbers, such 
that $\Cal H(b)<+\infty$, $\lim\limits_{n\rightarrow+\infty}b(n)=+\infty$ and 
$\lim\limits_{n\rightarrow+\infty}\frac{{\|T^n\|}_{L(X)}}{b(n)}=0$ (so that 
$r(T)\leq1$ by Theorem $6.1$). If $s:\Bbb N\rightarrow\Bbb R$ is any 
sequence of strictly positive real numbers, such that $\varDelta^ps$ is 
concave and unbounded from above for some $p\in\Bbb N$, and the 
sequence ${\bigl(\frac{b(n)}{s(n)}\bigr)}_{\!n\in\Bbb N}$ is bounded, then 
$\lim\limits_{n\rightarrow+\infty}\frac{{\|T^n\|}_{L(X)}}{s(n)}=0$, and each 
of conditions $(6.7.2)$--$(6.7.4)$ is equivalent to the following:
\roster
\item"(6.9.1)"the sequence ${\Biggl(\frac{\sum\limits_{k=0}^n(\varDelta s)(n-k)
\,T^k}{s(n)}\Biggr)}_{\!n\in\Bbb N}$ converges in $L(X)$.
\endroster
Finally, if $P\in L(X)$ is such that \ \,$\frac{\sum\limits_{k=0}^n(\varDelta s)
(n-k)\,T^k}{s(n)}\longrightarrow P$ in $L(X)$ as $n\rightarrow+\infty$ (so 
that conitions $(6.7.2)$--$(6.7.4)$ are also satisfied), then 
$P$ is the projection of $X$ onto $\Cal N(I_X-T)$ along $\Cal R(I_X-T)$.
\endproclaim
Let $\alpha\in(0,+\infty)$. By applying Theorem $6.7$ or Corollary $6.9$ 
to the sequences $b={\bigl((n+1)^\alpha\bigr)}_{n\in\Bbb N}$ and $s=
A_\alpha$ (see Example $4.9$, $(2.6)$ and Theorem $4.7$; see also 
$(2.5)$), we derive that if $T$ is a bounded linear operator on a complex 
nonzero Banach space $X$, such that $\lim\limits_{n\rightarrow+\infty}
\frac{{\|T^n\|}_{L(X)}}{n^\alpha}=0$, then the sequence ${\Biggl(\frac{\sum
\limits_{k=0}^nA_{\alpha-1}(n-k)\,T^k}{A_\alpha(n)}\Biggr)}_{\!n\in\Bbb
N}$ converges in $L(X)$ if and only if $1$ is either in $\rho(T)$, or a 
simple pole of $\goth R_T$. From this and from the result by 
E. Hille's mentioned in the Introduction ([Hi], Theorem $6$), together 
with Remark $6.8$, Theorem $2.9$ can be deduced. We remark that, 
however, Theorem $6.7$ does not completely extend Theorem $2.9$ to 
a larger class of sequences than the one of the sequences of Ces\`aro 
numbers (that is, the class of divergent nondecreasing sequences $s$ 
of strictly positive real numbers for which $\lim\limits_{n\rightarrow+\infty}
\frac{s(n+1)}{s(n)}=1$ and $\varDelta^qs\in\ell_1$ for some $q\in\Bbb N_2$), 
and neither does Corollary $6.9$ relative to the class of all sequences 
$s$ of strictly positive real numbers for which $\varDelta^ps$ is concave 
and unbounded from above for some $p\in\Bbb N$. Indeed, if $X$ is a complex 
nonzero Banach space, $T\in L(X)$, and $s$ is a nondecreasing 
sequence of strictly positive real numbers, such that $\lim\limits_
{n\rightarrow+\infty}s(n)=+\infty$, $\lim\limits_{n\rightarrow+\infty}\frac
{s(n+1)}{s(n)}=1$, $\varDelta^qs\in\ell_1$ for some $q\in\Bbb N_2$, 
and the sequence ${\Biggl(\frac{\sum\limits_{k=0}^n(\varDelta s)(n-k)\,
T^k}{s(n)}\Biggr)}_{\!n\in\Bbb N}$ converges in $L(X)$, then \,$\frac{
{\|T^n\|}_{L(X)}}{s(n)}$ need not converge to zero as $n\rightarrow+\infty$, 
even if $\varDelta^ps$ is concave and unbounded from above for some 
$p\in\Bbb N$. The following is an example.
\example{Example 6.10}
Let us consider the complex Banach space $\Bbb C^2$, endowed with 
the infinity norm (that is, ${\|(z_1,z_2)\|}_{\Bbb C^2}=\max\{|z_1|,|z_2|\}$ 
for all $(z_1,z_2)\in\Bbb C^2$). If $A\in L(\Bbb C^2)$ is the operator 
represented by the matrix $\pmatrix1&1\\0&1\endpmatrix$ (with 
respect to the canonical basis of $\Bbb C^2$), then $\sigma(A)=
\{1\}$. Furthermore, for each $n\in\Bbb N$, $A^n$ is represented by 
the matrix $\pmatrix1&n\\0&1\endpmatrix$.
\newline
Now let $T\in L(\Bbb C^2)$ be defined by $T=-A$. Then $\sigma(T)=
\{-1\}$, which gives $r(T)=1$ and $1\in\rho(T)$.
\newline
We define a sequence $a:\Bbb N\rightarrow\Bbb R$ of strictly positive 
real numbers as follows:
$$
a(n)=\cases1&\text{if }n=0\\\frac52&\text{if }n=1\\\frac1{n-1}+\frac2n+
\frac1{n+1}&\text{if }n\in\Bbb N_2.\endcases
$$
Since $a(2)=\frac73<\frac52=a(1)$, it follows that the sequence 
${\bigl(a(n))}_{n\in\Bbb Z_+}$ is strictly decreasing. Now let 
$s:\Bbb N\rightarrow\Bbb R$ be the sequence defined by $s=\varSigma
a$. Then $\varDelta s=a$. We also remark that $s(n)>0$ for each 
$n\in\Bbb N$, and $\lim\limits_{n\rightarrow+\infty}s(n)=+\infty$. 
Furthermore, since the sequence ${\bigl(s(n+1)-s(n)\bigr)}_{n\in\Bbb N}
={\bigl(a(n+1)\bigr)}_{n\in\Bbb N}$ is strictly decreasing, it follows that 
$s$ is concave. Then $s$ satisfies the hypotheses of Theorem $4.7$ 
for $p=0$ (so that $\lim\limits_{n\rightarrow+\infty}\frac{s(n+1)}{s(n)}=1$ 
and $\varDelta^2s\in\ell_1$).
\newline
We prove that
$$
\frac{\sum\limits_{k=0}^n(\varDelta s)(n-k)\,T^k}{s(n)}@>>n\rightarrow
+\infty>0_{L(\Bbb C^2)}\qquad\text{in \,}L(\Bbb C^2).\tag6.10.1
$$
We remark that, for each $k\in\Bbb N$, $T^k$ is represented by the 
matrix $\pmatrix{(-1)}^k&{(-1)}^kk\\0&{(-1)}^k\endpmatrix$. Hence 
proving $(6.10.1)$ is equivalent to proving that
$$
\frac{\sum\limits_{k=0}^n{(-1)}^k(\varDelta s)(n-k)}{s(n)}@>>n\rightarrow
+\infty>0\tag6.10.2
$$
and
$$
\frac{\sum\limits_{k=0}^n{(-1)}^kk\,(\varDelta s)(n-k)}{s(n)}@>>n\rightarrow
+\infty>0.\tag6.10.3
$$
We begin by proving $(6.10.2)$. We observe that for each $n\in\Bbb N$ 
we have
$$
\sum_{k=0}^n{(-1)}^k(\varDelta s)(n-k)=\sum_{k=0}^n{(-1)}^{n-k}
(\varDelta s)(k)={(-1)}^n\sum_{k=0}^n{(-1)}^ka(k).\tag6.10.4
$$
Since $a$ is eventually nonincreasing and $\lim\limits_{n\rightarrow+\infty}
a(n)=0$, we conclude that the series $\sum\limits_{n=0}^{+\infty}{(-1)}^na(n)$ 
converges, and consequently the sequence $\left(\,\sum\limits_
{k=0}^n{(-1)}^ka(k)\right)_{n\in\Bbb N}$ is bounded. Since $\lim\limits_
{n\rightarrow+\infty}s(n)=+\infty$, the desired result now follows from $(6.10.4)$.
\newline
Now we prove $(6.10.3)$. We first remark that, since for each $t\in
(-1,1)$ we have $\frac1{{(1-t)}^2}=\sum\limits_{n=1}^{+\infty}n\,t^{n-1}$, 
we obtain
$$
\sum_{n=0}^{+\infty}{(-1)}^nn\,t^n=-t\sum_{n=1}^{+\infty}n{(-t)}^{n-1}=
-\frac{t}{{(t+1)}^2}\qquad\text{for each }t\in(-1,1).\tag6.10.5
$$
Also, since $\sum\limits_{n=0}^{+\infty}\frac{t^n}{n+1}=\frac1t\sum\limits_
{n=1}^{+\infty}\frac{t^n}n=-\frac{\log(1-t)}t$ for each $t\in(0,1)$, we conclude 
that
$$
\multline
-\frac{{(t+1)}^2\log(1-t)}t=(t^2+2t+1)\sum_{n=0}^{+\infty}\frac{t^n}{n+1}
\\\allowdisplaybreak\hskip42pt
=\sum_{n=0}^{+\infty}\frac{t^{n+2}}{n+1}+\sum_{n=0}^{+\infty}\frac{2\,t^{n+1}
}{n+1}+\sum_{n=0}^{+\infty}\frac{t^n}{n+1}
=1+({\tsize\frac12}+2)t+\sum_{n=2}^{+\infty}\bigl({\tsize\frac1{n-1}+\frac2n
+\frac1{n+1}}\bigr)t^n
\\\allowdisplaybreak
=1+{\tsize\frac52}\,t+\sum_{n=2}^{+\infty}\bigl({\tsize\frac1{n-1}+\frac2n+
\frac1{n+1}}\bigr)t^n=\sum_{n=0}^{+\infty}a(n)\,t^n\qquad\text{for each }t\in
(0,1).
\endmultline
\tag6.10.6
$$
Since
$$
-\frac t{{(t+1)}^2}\left(-\frac{{(t+1)}^2\log(1-t)}t\right)=\log(1-t)=-\sum_{n=1}^
{+\infty}\frac{t^n}n\qquad\text{for each }t\in(0,1),
$$
from $(6.10.5)$ and $(6.10.6)$ it follows that
$$
\sum_{k=0}^n{(-1)}^kk\,a(n-k)=-\frac1n\qquad\text{for each }n\in\Bbb Z_+.
\tag6.10.7
$$
Since $\lim\limits_{n\rightarrow+\infty}s(n)=+\infty$, from $(6.10.7)$ 
we conclude that
$$
\frac{\sum\limits_{k=0}^n{(-1)}^kk\,(\varDelta s)(n-k)}{s(n)}=
\frac{\sum\limits_{k=0}^n{(-1)}^kk\,a(n-k)}{s(n)}=-\frac1{n\,s(n)}@>>
n\rightarrow+\infty>0,
$$
which is the desired result. We have thus proved $(6.10.2)$ and 
$(6.10.3)$, and consequently $(6.10.1)$. Hence $r(T)=1$, $1\in
\rho(T)$, and the sequence ${\Biggl(\frac{\sum\limits_{k=0}^n
(\varDelta s)(n-k)\,T^k}{s(n)}\Biggr)}_{\!n\in\Bbb N}$ converges in 
$L(\Bbb C^2)$ to $0_{L(\Bbb C^2)}$, which is, by the way, the projection 
of \,$\Bbb C^2$ onto \,$\{0_{\Bbb C^2}\}=\Cal N(I_{\Bbb C^2}-T)$ 
along \,$\Bbb C^2=\Cal R(I_{\Bbb C^2}-T)$. Nevertheless, we prove 
that the sequence ${\left(\frac{{\|T^n\|}_{L(\Bbb C^2)}}{s(n)}\right)}_
{n\in\Bbb N}$ does not converge to zero as $n\rightarrow+\infty$.
\newline
Since for each $n\in\Bbb N$ we have
$$
\multline
{\|T^n(z_1,z_2)\|}_{\Bbb C^2}={\|{(-1)}^n(z_1+nz_2,z_2)\|}_{\Bbb C^2}
=\max\{|z_1+nz_2|,|z_2|\}\\\allowdisplaybreak\leq\max
\{|z_1|+n|z_2|,|z_2|\}\leq(n+1)\max\{|z_1|,|z_2|\}=(n+1){\|(z_1,z_2)\|}_
{\Bbb C^2}
\endmultline
$$
and
$$
{\|T^n(1,1)\|}_{\Bbb C^2}={\|(n+1,1)\|}_{\Bbb C^2}=n+1,
$$
being ${\|(1,1)\|}_{\Bbb C^2}=1$ it follows that ${\|T^n\|}_{L(\Bbb C^2)}
=n+1$ for each $n\in\Bbb N$. On the other hand, since $s$ is 
concave, by Remark $3.4$ there exists $M\in(0,+\infty)$ such that \,
$\frac{s(n)}{n+1}\leq M$ for each $n\in\Bbb N$. Since $s(n)>0$ for 
each $n\in\Bbb N$, it follows that
$$
\frac{{\|T^n\|}_{L(\Bbb C^2)}}{s(n)}=\frac{n+1}{s(n)}\geq\frac1M
\qquad\text{for each }n\in\Bbb N.
$$
Hence the sequence ${\left(\frac{{\|T^n\|}_{L(\Bbb C^2)}}{s(n)}\right)}_
{n\in\Bbb N}$ does not converge to zero as $n\rightarrow+\infty$.
\newline
Actually, we can see that $\lim\limits_{n\rightarrow+\infty}\frac
{{\|T^n\|}_{L(\Bbb C^2)}}{s(n)}=+\infty$. Indeed, since for each 
$k\in\Bbb N_2$ we have $\frac1k\leq\frac1x$ for each $x\in[k-1,k]$, 
and consequently $\frac1k\leq\int_{k-1}^k\frac1x\,dx$, it follows 
that for each $n\in\Bbb N_3$ we have
$$
\multline
s(n)=1+{\tsize\frac52}+\sum_{k=2}^n\bigl({\tsize\frac1{k-1}+\frac2
k+\frac1{k+1}}\bigr)\leq{\tsize\frac72}+4\sum_{k=2}^n{\tsize\frac1
{k-1}}={\tsize\frac72}+4+4\sum_{k=3}^n{\tsize\frac1{k-1}}={\tsize
\frac{15}2}+4\sum_{k=2}^{n-1}{\tsize\frac1k}\\\allowdisplaybreak
\leq{\tsize\frac{15}2}+4\sum_{k=2}^{n-1}\,\int_{k-1}^k{\tsize\frac1x}
\,dx={\tsize\frac{15}2}+4\int_1^{n-1}{\tsize\frac1x}\,dx={\tsize\frac
{15}2}+4\log(n-1).
\endmultline
$$
Hence
$$
\frac{{\|T^n\|}_{L(\Bbb C^2)}}{s(n)}=\frac{n+1}{s(n)}\geq\frac{n+1}
{\frac{15}2+4\log(n-1)}\qquad\text{for each }n\in\Bbb N_3,
$$
which gives $\lim\limits_{n\rightarrow+\infty}\frac{{\|T^n\|}_
{L(\Bbb C^2)}}{s(n)}=+\infty$.
\endexample

\enddocument